\documentclass[preprint,11pt]{elsarticle}

\usepackage{ifpdf}
\usepackage{amsmath}
\usepackage{amsfonts}
\usepackage{amsthm}
\usepackage{mathrsfs}
\usepackage{color}
\usepackage{graphicx}
\usepackage{subfigure}
\usepackage{float}
\usepackage{overpic}
\usepackage{threeparttable}
\usepackage{epstopdf}
\usepackage{dcolumn}
\usepackage{multirow}
\usepackage{booktabs}
\biboptions{numbers,sort&compress}

\usepackage{graphicx,natbib,amssymb,lineno}
\ifpdf
\usepackage[%
  pdftitle={Instructions for use of the document class
    elsart},%
  pdfauthor={},%
  pdfsubject={The preprint document class elsart},%
  pdfkeywords={instructions for use, elsart, document class},%
  pdfstartview=FitH,%
  bookmarks=true,%
  bookmarksopen=true,%
  breaklinks=true,%
  colorlinks=true,%
  linkcolor=blue,anchorcolor=blue,%
  citecolor=blue,filecolor=blue,%
  menucolor=blue,pagecolor=blue,%
  urlcolor=blue]{hyperref}
\else
\usepackage[%
  breaklinks=true,%
  colorlinks=true,%
  linkcolor=blue,anchorcolor=blue,%
  citecolor=blue,filecolor=blue,%
  menucolor=blue,pagecolor=blue,%
  urlcolor=blue]{hyperref}
\fi

\makeatletter
\def\elsartstyle{%
    \def\normalsize{\@setfontsize\normalsize\@xiipt{14.5}}
    \def\small{\@setfontsize\small\@xipt{13.6}}
    \let\footnotesize=\small
    \def\large{\@setfontsize\large\@xivpt{18}}
    \def\Large{\@setfontsize\Large\@xviipt{22}}
    \skip\@mpfootins = 18\p@ \@plus 2\p@
    \normalsize
} \@ifundefined{square}{}{} \makeatother

\newtheorem{theorem}{Theorem}[section]
\newtheorem{lemma}[theorem]{Lemma}
\newtheorem{proposition}[theorem]{Proposition}

\theoremstyle{definition}

\newtheorem{example}[theorem]{Example}

\theoremstyle{remark}

\makeatletter
\def\ps@pprintTitle{%
  \let\@oddhead\@empty
  \let\@evenhead\@empty
  \let\@oddfoot\@empty
  \let\@evenfoot\@oddfoot
}
\makeatother

\def\ps@pprintTitle{%
  \let\@oddhead\@empty
  \let\@evenhead\@empty
  \def\@oddfoot{\reset@font\hfil\thepage\hfil}
  \let\@evenfoot\@oddfoot
}
\begin{document}

\begin{frontmatter}

\title{Nilpotent center conditions in cubic switching polynomial Li\'enard systems by higher-order analysis}

\author[myaddress,mymainaddress]{Ting Chen}
\ead{chenting0715@126.com}
\author[mysecondaryaddress]{Feng Li}
\ead{lf0539@126.com}
\author[mythirdaddress]{Pei Yu\corref{mycorrespondingauthor}}
\cortext[mycorrespondingauthor]{Corresponding author}
\ead{pyu@uwo.ca}
\address[myaddress]{College of Science, National University of Defense Technology, Changsha, 410073,  China  \vspace*{0.05in}}
\address[mymainaddress]{School of Statistics and Mathematics, Guangdong
University of Finance and Economics,\\
Guangzhou, 510320, China \vspace*{0.05in} }
\address[mysecondaryaddress]{School of Mathematics and Statistics,
Linyi University, Linyi, 276005, China \vspace*{0.05in}}
\address[mythirdaddress]{Department of Mathematics,
Western University, London, Ontario, N6A 5B7, Canada}

\begin{abstract}
The aim of this paper is to investigate two classical problems related to
nilpotent center conditions and bifurcation of limit cycles
in switching polynomial systems.
Due to the difficulty in calculating the Lyapunov constants of
switching polynomial systems at non-elementary singular points,
it is extremely difficult to use the existing Poincar\'e-Lyapunov method
to study these two problems.
In this paper, we develop a higher-order Poincar\'e-Lyapunov method to
consider the nilpotent center problem in switching polynomial systems,
with particular attention focused on cubic switching Li\'enard systems.
With proper perturbations, explicit center conditions are derived
for switching Li\'enard systems at a nilpotent center,
which is characterized as global.
Moreover, with Bogdanov-Takens bifurcation theory, the existence of
five limit cycles around the nilpotent center is proved
for a class of switching Li\'enard systems, which is a
new lower bound of cyclicity for such polynomial systems
around a nilpotent center.
\end{abstract}


\end{frontmatter}

\section{Introduction}

During the past several decades, a large number of works has been
focused on the study of the so-called Li\'enard equation,
\begin{equation}\label{Eqn1-1}
\ddot{x}+f(x)\dot{x}+g(x)=0,
\end{equation}
which frequently appears in many disciplines and applications
\cite{Gasull1989}. The dot in \eqref{Eqn1-1}
denotes derivative with respect to time $t$.
By introducing $y=F(x)+\dot{x}=\int_0^x f(x) dx+\dot{x}$,
the equation \eqref{Eqn1-1} can be brought to the following planar
first-order differential system,
\begin{equation}\label{Eqn1-2}
\begin{array}{ll}
\dot{x}= y-F(x),\\
\dot{y}= -\,g(x),
\end{array}
\end{equation}
which is the so-called generalized Li\'enard system.

For the Li\'enard system \eqref{Eqn1-2}, three important problems related
to qualitative dynamical behaviours are classified as the center conditions,
the number of limit cycles and the global phase portraits,
see \cite{CHB2015,Cherkas1977,LlibreV2022,ZDHD1992}.
Christopher \cite{Christopher} introduced an algebraic approach
to classify linear-type centers in smooth polynomial Li\'enard systems,
which are called elementary singular points characterized by a pair of purely
imaginary eigenvalues. Further, such a center is considered as global
if the whole vector field of the system is filled with periodic orbits
except for this point. Later, Llibre and Valls \cite{LlibreV2022}
studied all types of generalized Li\'enard systems
having a global linear-type center.

Assume that the polynomials $F(x)$ and $g(x)$ are given in
the form of $F(x)=\sum^n_{i=0}a_ix^i$ and $g(x)=\sum^m_{i=0}b_ix^i$.
The problem of the number of bifurcating limit cycles and their relative
positions for the classical Li\'enard system (i.e., $g(x)=x$)
is the well-known Smale's 13th problem \cite{Smale,Caubergh}.
The authors of \cite{LMP1977} conjectured that the classical Li\'enard system
can have at most $[\frac{n-1}{2}]$ limit cycles, where $[\cdot]$ denotes
the integer function.
Maesschalck and Dumortier \cite{Maesschalck2011} proved that some
classical Li\'enard system can have at least $[\frac{n-1}{2}]+2$ limit
cycles for $n\geq6$.

It has been noted that much attention for the generalized Li\'enard system
was paid to consider the maximum number of limit cycles
bifurcating from a monodromic singular point, which is classified as
either a center or a focus, usually called Hopf cyclicity,
see for example \cite{ChL1999,TH2011,JHYL2007,JH2009}.
For the Hopf cyclicity at the origin of \eqref{Eqn1-2},
Han \cite{Han2013} proved it to be $[\frac{n-1}{2}]+2$
when $g(-x)=-g(x)$; Tian and Han \cite{TH2011} showed it
to be $[\frac{2n-1}{3}]$ when $\text{deg}(g)=2$;
and Tian {\it et al.} \cite{THX2019} proved it to be
$[\frac{3n-1}{4}]$ when $\text{deg}(g)=3$.
Chen {\it el at.} \cite{CHB2016,CHB2020} proved
the existence of two limit cycles in two classes of cubic Li\'enard systems,
and analyzed their global dynamics.
However, to the best of our knowledge, the three problems
for the general Li\'enard system are still open.

In recent years, increasing attention has been attracted to
the research in dynamics of non-smooth systems since
non-smoothness has been included more and more in models describing
problems arising in engineering \cite{CDT2018}, epidemiology \cite{Tang}
and electronics \cite{Colombo}, and in particular in
the generalized Li\'enard system \cite{Adriana,LT2015,MM2014,SHR2016}.
For the non-smooth Li\'enard system \eqref{Eqn1-2},
$F(x)$ and $g(x)$ are usually assumed to be piecewise smooth.
For example, Chen {\it el at.} \cite{CDT2018} studied the global dynamics
of a mechanical system with dry friction, which can be transformed to a
piecewise smooth Li\'enard system.
However, the center and Hopf cyclicity problems become
extremely complicated in switching systems.
To overcome the difficulty, the authors of \cite{Gasull} developed
a useful approach for computing the Lyapunov constants for the
switching polynomial systems with an elementary singular point.
With this method, they solved the linear-type center problem
for a class of switching Li\'enard systems.
In \cite{TH2011,THX2019}, the Hopf cyclicities are obtained
for system \eqref{Eqn1-2} when $F(x)$ is a piecewise polynomial
having a switching manifold at $x = 0$ with $\text{deg}(g)=2,\,3$.

It should be pointed out that although many research results
studies have been obtained
on the center problem and bifurcation of small-amplitude limit cycles
for switching Li\'enard systems associated with elementary singular points,
no attention has been paid to switching Li\'enard systems
associated with an isolated nilpotent singular point.
By an isolated nilpotent singular point in planar polynomial systems,
it means that the two eigenvalues of the Jacobian matrix of the system,
evaluated at the singular point,
are zero but the Jacobian matrix is not null, more details can be found in
\cite{Garcia,Giacomini,F.Li2,F.Li2021,Liu2014,Liu2011,Strozyna,Yang,YL2017}.
From Theorem 3.5 in \cite{Dumortier2006} we know that if $m$ is the
smallest integer satisfying $b_m\neq0$, the multiplicity of the nilpotent
singular point of system \eqref{Eqn1-2} is exactly $m$.
Further, let $n$ be the smallest number for which $a_n\neq0$.
Then, the local qualitative properties on the nilpotent origin $(0,0)$
of \eqref{Eqn1-2} are summarized
in Table \ref{T1}.

\begin{table}[!h]\label{T1}
\caption{The local qualitative properties on the nilpotent singular point
$(0,0)$ of \eqref{Eqn1-2}.}
{\scriptsize\begin{center}
{\begin{tabular}{|l|l|l|l|l|}
\hline
\multicolumn{4}{|c|}{Conditions}  &Type of the origin  \\
\hline
\multirow{3}*{$a_{n}=0$, $b_m\neq0$} &\multicolumn{3}{c|}{$m$ odd,
$b_m>0$} &a center or a focus \\[1mm]
\cline{2-5}
&\multicolumn{3}{c|}{$m$ odd, $b_m<0$} & a saddle \\[1mm]
\cline{2-5}
&\multicolumn{3}{c|}{$m$ even}  & a cusp\\[1mm]
\hline
\multirow{8}*{$\begin{array}{l} \!\!\! a_{n}\neq0, \ b_m\neq0, \\[2.0ex]
\end{array}$} &\multirow{6}*{$m=2k \!+\!1$}
& \multicolumn{2}{c|}{$b_m<0$} &a saddle\\[1.0mm]
\cline{3-5}
\multirow{7}*{$\Delta=-4(n \!+\! 1)b_{m} \!+\! n^2a_n^2$} &
&\multirow{5}*{$b_m>0$}
& \multirow{3}*{$k>n-1$, or $k=n-1$ and $\Delta\geq0$, $n$ even}
& consisting of one \\[-0.5mm]
& & & & hyperbolic sector and \\[-0.5mm]
& & & & one elliptic sector \\[0.5mm]
\cline{4-5}
& & &$k>n-1$, or $k=n-1$ and $\Delta\geq0$, $n$ odd & a node \\[1mm]
\cline{4-5}
& & &$k<n-1$, or $k=n-1$ and $\Delta<0$ & a center or a focus\\[1mm]
\cline{2-5}
&\multirow{2}*{$m=2k$} &\multicolumn{2}{c|}{$n<k$}  & a saddle-node\\[1mm]
\cline{3-5}
& &  \multicolumn{2}{c|}{$n\geq k$}  & a cusp\\[1mm]
\hline
\end{tabular}}
\end{center}}
\end{table}

Regarding Table \ref{T1},
it is worth to mention that the origin of \eqref{Eqn1-2}
is a monodromic singular point and a cusp when the parameters satisfy
\begin{equation}\label{X1}
\begin{aligned}
\Xi_1=&\ \big\{(a_n,b_m)\in\mathbb{R}^2:\ a_n=0,\ b_{2k+1}>0\big\}\bigcup
\big\{(a_n,b_m)\in\mathbb{R}^2:\ a_n\neq0,\ b_{2k+1}>0,\ k>n\big\}\\
&\ \bigcup
\big\{(a_n,b_m)\in\mathbb{R}^2:\ a_n\neq0,\ b_{2k+1}>0,\ k=n,\ \Delta<0\big\}
\end{aligned}
\end{equation}and
\begin{equation}\label{X2}
\begin{aligned}
\Xi_2=\big\{(a_n,b_m)\in\mathbb{R}^2:\,\,a_n=0,\,\,b_{2k}\neq0\big\}\bigcup
\big\{(a_n,b_m)\in\mathbb{R}^2:\,\,a_n\neq0,\,\,b_{2k}\neq0,\,\,k\leq n\big\},
\end{aligned}
\end{equation}respectively.

In this paper, we will develop a higher-order Poincar\'e-Lyapunov method
to determine the nilpotent center conditions and bifurcation of
small-amplitude limit cycles in switching polynomial systems.
It is a more challenging and interesting compared to the works for switching
systems with elementary singular points. We will apply this method to investigate
a class of cubic switching polynomial Li\'enard systems with a nilpotent
singular point. Without loss of generality, the cubic switching
Li\'enard systems can be written in the form of the differential equations,
\begin{equation}\label{Eqn1-3}
\left(\begin{array}{cc}
\dot{x}\\
\dot{y}
\end{array}\right)
=\left\{\begin{aligned}&\left(\begin{array}{c}
y-(a_0^++a_1^+x+a_2^+x^2+a_3^+x^3) \\[0.5ex]
-(b_0^++b_1^+x+b_2^+x^2+b_3^+x^3)
\end{array}\right), &\text{if} \ \ x\geq0,\\
&\left(\begin{array}{c}
y-(a_0^-+a_1^-x+a_2^-x^2+a_3^-x^3) \\[0.5ex]
-(b_0^-+b_1^-x+b_2^-x^2+b_3^-x^3)
\end{array}\right), &\text{if} \ \ x<0,
\end{aligned}
\right.
\end{equation}
where $x=0$ is the unique switching manifold, and
$\lambda \!=\! (a_i^{\pm},b_j^{\pm})
\! \in \! \mathbf{R}^{16}$, $i=0,1,2,3$, represents
the parameter vector.
For the convenience in the following analysis, we call the system
with ``$+$'' sign ``the first system''
and the system with ``$-$'' sign ``the second system''.

Assume that $(0,0)$ is a singular point of \eqref{Eqn1-3}.
Then, we have $a_0^{\pm}=b_0^{\pm}=0$.
Thus, the Jacobian matrices of the first and the second systems
of \eqref{Eqn1-3} evaluated at the origin are given by
\begin{equation}\label{Eqn1-4}
J^{\pm}=\left[\begin{array}{cc}
-a_1^{\pm}& 1 \\
-b_1^{\pm}& 0
\end{array}\right].
\end{equation}
The necessary and sufficient conditions
for the origin of \eqref{Eqn1-3} to be an isolated nilpotent singular point
are ${\rm Tr}(J^{\pm})=\det (J^{\pm})=0$, with
$J^{\pm}$ not being identically zero.
It is easy to obtain that $a_1^{\pm}=b_1^{\pm}=0$.
Then, \eqref{Eqn1-3} is reduced to
\begin{equation}\label{Eqn1-5}
\left(\begin{array}{cc}
\dot{x}\\
\dot{y}
\end{array}\right)
=\left\{\begin{aligned}&\left(\begin{array}{c}
y-(a_2^+x^2+a_3^+x^3) \\[0.5ex]
-(b_2^+x^2+b_3^+x^3)
\end{array}\right), &\text{if}~~x\geq0,\\
&\left(\begin{array}{c}
y-(a_2^-x^2+a_3^-x^3) \\[0.5ex]
-(b_2^-x^2+b_3^-x^3)
\end{array}\right), &\text{if}~~x<0.
\end{aligned}
\right.
\end{equation}

It follows from Table \ref{T1} that different types of the nilpotent singular
point can generate much more rich dynamics than that of the elementary one,
leading to that the determination of the center conditions of
\eqref{Eqn1-5} becomes more involved. In fact, the center of system
\eqref{Eqn1-5} on the switching manifold
can be classified as two monodromic singular points, or a cusp
and a monodromic singular point, or two cusps.
Due to this complexity it is extremely difficult to consider
the nilpotent center problem in the switching polynomial Li\'enard systems
of general degree $n$.
Here, we derive the necessary and sufficient conditions for the
center problem associated with the nilpotent origin of the
cubic switching polynomial Li\'enard system \eqref{Eqn1-5}.
We have the following result.

\begin{theorem}\label{Th1}
Assume $b_2^+\geq0$ and $b_2^-\leq 0$.
The origin of the cubic switching Li\'enard system \eqref{Eqn1-5} is a nilpotent center if one of the following conditions holds:
\begin{equation}\label{Eqn1-6}
\begin{aligned}
\mathrm{I}:&\ \ a_2^{\pm}=a_3^{\pm}=b_2^-=0,\ b_2^+>0,\ b_3^->0;\\
\mathrm{II}:&\ \ a_2^{\pm}=a_3^{\pm}=0,\ b_2^+>0,\ b_2^-<0;\\
\mathrm{III}:&\ \ a_2^{\pm}=b_3^{\pm}=a_3^+b_2^--a_3^-b_2^+=0,\
a_3^+(b_2^++b_2^-)\neq0 ,\ b_2^+>0,\ b_2^-<0;\\
\mathrm{IV}:&\ \ a_2^-=a_2^+,\ a_3^-=-a_3^+,\ b_2^-=-b_2^+<0,\ b_3^-=b_3^+;\\
\mathrm{V}:&\ \ a_2^{\pm}=a_3^{\pm}=b_2^+=0,\ b_2^-<0,\ b_3^+>0;\\
\mathrm{VI}:&\ \ a_2^{\pm}=b_2^{\pm}=0,\ a_3^-=-a_3^+,\
b_3^-=b_3^+> \tfrac{3}{4} (a_3^+)^2;\\
\mathrm{VII}:&\ \ a_2^{\pm}=a_3^{\pm}=b_2^{\pm}=0,\ b_3^+>0,\ b_3^->0;\\
\mathrm{VIII}:&\ \ a_2^-=a_2^+\neq0,\ a_3^-=-a_3^+,\ b_2^{\pm}=0,\
b_3^-=b_3^+> \tfrac{1}{3} (a_2^+)^2.\\
\end{aligned}
\end{equation}
\end{theorem}

The next result further characterizes the nilpotent center $(0,0)$ of
the cubic switching Li\'enard system \eqref{Eqn1-5} to be global.

\begin{theorem}\label{Th2}
Assume $b_2^+\geq0$ and $b_2^-\leq 0$.
The origin of the cubic switching Li\'enard system \eqref{Eqn1-5} is a nilpotent global center
if one of the following conditions holds:
\begin{equation}\label{Eqn1-7}
\begin{aligned}
\mathrm{G_1}:&\ \ a_2^{\pm}=a_3^{\pm}=0,\ b_3^+>0,\ b_3^->0;\\
\mathrm{G_2}:&\ \ a_2^-=a_2^+,\ a_3^{\pm}=0,\ b_2^-=-b_2^+,\ b_3^-=b_3^+>\tfrac{1}{2}(a_2^+)^2.\\
\end{aligned}
\end{equation}
\end{theorem}

Moreover, regarding the bifurcation of limit cycles around the
origin of \eqref{Eqn1-5}, we construct a perturbed system
using the center condition $\mathrm{IV}$
and obtain that the maximal number of bifurcating limit cycles
from the nilpotent origin of system \eqref{Eqn1-5} is $5$.
This is a new lower bound on the limit cycles for such cubic
switching Li\'enard systems around a nilpotent singular point.

\begin{theorem}\label{Th3}
With the nilpotent center condition $\mathrm{IV}$ in Theorem {\rm \ref{Th1}},
the cubic switching Li\'enard system \eqref{Eqn1-5} can have at least $5$ limit cycles bifurcating from the origin
by higher-order cubic perturbations.
\end{theorem}

The rest of the paper is organized as follows.
In the next section, we present our high-order Poincar\'{e}-Lyapunov method
with some formulas which are needed in Section 3.
Section 3 is devoted to derive the nilpotent center conditions
at the origin of system \eqref{Eqn1-5}. The conditions
on the global nilpotent center of system \eqref{Eqn1-5} are obtained
in Section 4. In Section 5, a perturbed system of \eqref{Eqn1-5} is
constructed to show the bifurcation of $5$ limit cycles
from the origin of \eqref{Eqn1-5}.
Finally, conclusion is drawn in Section 6.

\section{The high-order Poincar\'e-Lyapunov method}

We consider the following switching nilpotent systems divided by the $y$-axis,
{\begin{equation}\label{Eqn2-1}
(\dot{x},\, \dot{y})=\left\{\begin{aligned}
\Big(y+\sum\limits_{i+j=2}^{n}A_{ij}^+x^iy^j, \
\sum\limits_{i+j=2}^{n}B_{ij}^+x^iy^j\Big), \ \ \text{if} \ \
x\geq0,\\[0.0ex]
\Big(y+\sum\limits_{i+j=2}^{n}A_{ij}^-x^iy^i, \
\sum\limits_{i+j=2}^{n}B_{ij}^-x^iy^j\Big), \ \ \text{if} \ \
x<0,
\end{aligned}\right.
\end{equation}where $A_{ij}^{\pm}$ and $B_{ij}^{\pm}$ are parameters.
We have the following proposition
for proving the nilpotent center conditions at the origin $(0,0)$, which is the common nilpotent singular point in both the first and
the second systems of \eqref{Eqn2-1}.

\begin{proposition}\label{p11}
Assume that the origin $(0,0)$ of the switching nilpotent system \eqref{Eqn2-1} is monodromic.
If the first and second systems in \eqref{Eqn2-1} have respectively
the first integrals $I^+(x,y)$ and $I^-(x,y)$ near the origin, and
either both $I^+(0,y)$ and $I^-(0,y)$ are even functions in $x$ or
$I^+(0,y)\equiv I^-(0,y)$, then $(0,0)$  is
a nilpotent center.
\end{proposition}

See \cite{Gasull} for more details about Proposition \ref{p11}.
In \cite{LY2015}, the authors redefined the symmetry of switching systems.
By modifying the conditions, we obtain the following result for
proving the nilpotent center conditions of \eqref{Eqn2-1} at the origin.

\begin{proposition}\label{p22}
Assume that the nilpotent origin of the switching nilpotent system \eqref{Eqn2-1} is monodromic.
If the systems \eqref{Eqn2-1} are symmetric
with respect to the $x$-axis, i.e., the parameters on the right-hand side
of \eqref{Eqn2-1} satisfy
\begin{equation}\label{Eqn2-2}
\begin{aligned}
A_{i,2k}^{\pm}=B_{i,2k+1}^{\pm}=0,
\end{aligned}
\end{equation}
or the systems in \eqref{Eqn2-1} are symmetric with respect to the $y$-axis,
i.e., the parameters on the right-hand side of \eqref{Eqn2-1} satisfy
\begin{equation}\label{Eqn2-3}
\begin{aligned}
A_{2k+1,j}^+=-A_{2k+1,j}^-, \quad A_{2k,j}^+=A_{2k,j}^-, \quad B_{2k+1,j}^+=B_{2k+1,j}^-, \quad B_{2k,j}^+=-B_{2k,j}^-,
\end{aligned}
\end{equation}
then the origin of \eqref{Eqn2-1} is a nilpotent center.
\end{proposition}

Now we introduce our higher-order Poincar\'{e}-Lyapunov method to
study the switching system \eqref{Eqn2-1}, which establishes the
relation between the unperturbed systems and the perturbed systems based on
Bogdanov-Takens (B-T) bifurcation theory.
Hence, we consider the following perturbed system of \eqref{Eqn2-1},
{\begin{equation}\label{Eqn2-4}
\left(\begin{array}{cc}
\dot{x}\\
\dot{y}
\end{array}\right)
=\left\{\begin{aligned}
&\left(\begin{array}{c}
\vartheta_1 \varepsilon x+y+\displaystyle\sum\limits_{i+j=2}^{n}A_{ij}^+x^iy^j
+\displaystyle\sum\limits_{k=1,i+j=2}^{n}P_{kij}^+\varepsilon^k x^iy^j \\[2.0ex]
-\varepsilon^2x+\vartheta_2\varepsilon y
+\displaystyle\sum\limits_{i+j=2}^{n}B_{ij}^+x^iy^j
+\displaystyle\sum\limits_{k=1,i+j=2}^{n}Q_{kij}^+\varepsilon^kx^iy^j
\end{array}\right), \ \ &\text{if} \ \ x\geq0,\\[1.0ex]
&\left(\begin{array}{c}
\vartheta_1 \varepsilon x+y+\displaystyle\sum\limits_{i+j=2}^{n}A_{ij}^-x^iy^i
+\displaystyle\sum\limits_{k=1,i+j=2}^{n}P_{kij}^-\varepsilon^kx^iy^j \\[2.0ex]
-\varepsilon^2x+\vartheta_2\varepsilon y
+\displaystyle\sum\limits_{i+j=2}^{n}B_{ij}^-x^iy^j
+\displaystyle\sum\limits_{k=1,i+j=2}^{n}Q_{kij}^-\varepsilon^kx^iy^j
\end{array}\right), \ \ &\text{if} \ \ x<0,
\end{aligned}\right.
\end{equation}
where $-\varepsilon^2 x$ is called unfolding with sufficiently
small $\varepsilon>0$, $\lambda^{\pm}=(\vartheta_{1,2},A_{ij}^{\pm},
B_{ij}^{\pm},P_{ij}^{\pm},Q_{ij}^{\pm})$ represent two parameter vectors.
For convenience, we denote that \eqref{Eqn2-1} is the {\it limit} system
of \eqref{Eqn2-4}.

Note that the linear perturbation terms used in \eqref{Eqn2-4} are
$\varepsilon^2 x$ rather than $\varepsilon x$ can avoid the
$\sqrt{\varepsilon}$ and $\frac{1}{\varepsilon^{k}}$ perturbation terms ($k\in \mathbb{Q^+}$) in later
transformed systems. Based on the
relation established for these two systems \eqref{Eqn2-1} and \eqref{Eqn2-4},
and the B-T bifurcation theory,
we directly have the following lemma. More detailed discussions on
this subject can be found in \cite{Chen4}.

\begin{lemma}\label{L0}
Assume that the origin of \eqref{Eqn2-1} is monodromic, and that
the two systems in \eqref{Eqn2-1} are limit systems of \eqref{Eqn2-4}.
If the linear-type elementary singular point at the origin of \eqref{Eqn2-4}
becomes a linear-type center, then the origin of \eqref{Eqn2-1} is a nilpotent center.
\end{lemma}

We give the following example to illustrate our basic idea,
as it is known that the nilpotent polynomial smooth systems
can be transformed to the normal form (see \cite{HY2012}),
\begin{equation}\label{Eqn2-5}
\begin{aligned}
\dot{x}=&\ y,\\
\dot{y}=&\ a_k \big( 1+h(x) \big)x^k+b_ky \big(1+g(x) \big)x^{k-1}+y^2p(x,y),
\end{aligned}
\end{equation}
where $k\geq2$, $h(x)$, $g(x)$ and $p(x,y)$ are real analytic functions
satisfying $h(0)=g(0)=0$. For example, we consider a codimension-2
symmetric B-T bifurcation of the cubic normal form \eqref{Eqn2-5},
given in the following form:
\begin{equation}\label{Eqn2-6}
\begin{aligned}
\dot{x}=&\ y,\\
\dot{y}=&\ \varepsilon_1 x+\varepsilon_2y+a_3x^3+b_3x^2y,
\end{aligned}
\end{equation}
where the terms $\varepsilon_1 x+\varepsilon_2y$ are called unfolding
with small $\varepsilon_1$ and $\varepsilon_2$.
Note that if $\varepsilon_1<0$, the system \eqref{Eqn2-6} can have
Hopf bifurcation near the origin from the bifurcation line $\varepsilon_2=0$.
It is easy to check that the elementary origin of \eqref{Eqn2-6}
is a center if and only if $b_{3}=0$. The isolated elementary origin
of \eqref{Eqn2-6} is reduced to a nilpotent point when
$\varepsilon_1\rightarrow0^{-}$ and $\varepsilon_2\rightarrow0$.
Furthermore, it can be verified that the nilpotent monodromic origin
(when $a_3<0$) of system \eqref{Eqn2-6} without unfolding
is also a center if $b_3=0$.

The difficulty arising from the problem of distinguishing a center from
a focus in switching nilpotent system \eqref{Eqn2-1} is that the problem
may be not algebraically solvable. That is, it does not have an infinite
sequence of independent polynomials involving the coefficients of the
systems such that the Lyapunov constants vanish simultaneously,
which guarantees the existence of a center.
We will show that the higher-order Poincar\'e-Lyapunov method
we develop for the switching polynomial system \eqref{Eqn2-4} to
determine nilpotent center can overcome the difficulty.
More details are described below.

To achieve this, introducing the transformation
$(x,y,t)\rightarrow (\varepsilon^3 x, \varepsilon^2 y, \frac{t}{\varepsilon})$
into \eqref{Eqn2-4}, we obtain
{\begin{equation}\label{Eqn2-7}
\left(\begin{array}{cc}
\dot{x}\\
\dot{y}
\end{array}\right)
=\left\{\begin{aligned}
&\left(\begin{array}{c}
\vartheta_2x-y+\displaystyle\sum\limits_{k=1,i+j=2}^{n}
\widetilde{A_{kij}}^+(\lambda^+)\varepsilon^kx^iy^j  \\[2.5ex]
x+\vartheta_1y+\displaystyle\sum\limits_{k-1,i+j=2}^{n}
\widetilde{B_{kij}}^+(\lambda^+)\varepsilon^kx^iy^j
\end{array}\right), \ \ &\text{if} \ \ y\geq0,\\[0.5ex]
&\left(\begin{array}{c}
\vartheta_2x-y+\displaystyle\sum\limits_{k=1,i+j=2}^{n}
\widetilde{A_{kij}}^-(\lambda^-)\varepsilon^kx^iy^j \\[2.5ex]
x+\vartheta_1y+\displaystyle\sum\limits_{k=1,i+j=2}^{n}
\widetilde{B_{kij}}^-(\lambda^-)\varepsilon^kx^iy^j
\end{array}\right), \ \ &\text{if} \ \ y<0,
\end{aligned}\right.
\end{equation}
where $\widetilde{A_{kij}}^{\pm}(\lambda^{\pm})$ and
$\widetilde{B_{kij}}^{\pm}(\lambda^{\pm})$ are the polynomials
in the parameter vector $\lambda^{\pm}$.

With the polar coordinates transformation:
$x=r \cos\theta$ and $y=r \sin\theta$, the perturbation system \eqref{Eqn2-7}
becomes the form,
\begin{equation}\label{Eqn2-8}
\frac{{\rm d}r}{{\rm d}\theta}=\left\{\begin{aligned}
&\frac{\displaystyle\sum\limits_{k=0,i+j=1}^nc_{kij}^+(\lambda^+)
\cos\theta^i\sin\theta^j\varepsilon^kr^{i+j}}
{1+\displaystyle\sum\limits_{k=0,i+j=1}^nd_{kij}^+(\lambda^+)\cos\theta^i
\sin\theta^j\varepsilon^kr^{i+j-1}}, \quad
\text{if} \ \ \theta\in[0,\pi],\\[1.0ex]
&\frac{\displaystyle\sum\limits_{k=0,i+j=1}^{n}c_{kij}^-(\lambda^-)
\cos\theta^i\sin\theta^j\varepsilon^kr^{i+j}}
{1+\displaystyle\sum\limits_{k=0,i+j=1}^{n}d_{ij}^-(\lambda^-)\cos\theta^i
\sin\theta^j\varepsilon^kr^{i+j-1}}, \quad \text{if} \ \ \theta\in(\pi,2\pi),
\end{aligned} \right.
\end{equation}
where $c_{kij}^{\pm}(\lambda^{\pm})$ and
$d_{kij}^{\pm}(\lambda^{\pm})$ are the polynomials in
the parameter vector $\lambda^{\pm}$.

Let $r^+(\xi,\lambda^+,\varepsilon,\theta)
=\sum_{k\geq1}v_{k}^+(\lambda^+,\varepsilon,\theta)\xi^k$ and
$r^-(\xi,\lambda^-,\varepsilon,\theta)
=\sum_{k\geq1}v_{k}^-(\lambda^+,\varepsilon,\theta)\xi^k$
be the solutions of the first and second systems of \eqref{Eqn2-8}
associated with the initial conditions
$r^+(\xi,\lambda^+,\varepsilon,0)=r^-(\xi,\lambda^-,\varepsilon,\pi)=\xi$.
We denote by
$$
\Pi^+(\xi)=r^+(\xi,\lambda^+,\varepsilon,\pi)
=\sum\limits_{k\geq1}v_{k}^+(\lambda^+,\varepsilon,\pi)\xi^k
$$
and
$$
\Pi^-(\xi)=r^-(\xi,\lambda^+,\varepsilon,2\pi)
=\sum\limits_{k\geq1}v_{k}^-(\lambda^-,\varepsilon,\pi)\xi^k,
$$
the first half-return map $\Pi^+(\xi)$ and the second half-return
map $\Pi^-(\xi)$, respectively,
where $v_k^{\pm}(\lambda^{\pm},\varepsilon,\pi)$ are the coefficients
in Taylor expansions. However, it is extremely difficult to composite these
two maps to compute the displacement map of \eqref{Eqn2-8}.
We may follow the procedure in \cite{Gasull} and introduce
the transformation $(x,y,t)\rightarrow(x,-y,-t)$ into
the piecewise smooth system \eqref{Eqn2-7} to yield
{\begin{equation}\label{Eqn2-7-0}
\left(\begin{array}{cc}
\dot{x}\\
\dot{y}
\end{array}\right)
=\left\{\begin{aligned}
&\left(\begin{array}{c}
-\vartheta_2x-y-\displaystyle\sum\limits_{k=1,i+j=2}^{n}
\widetilde{A_{kij}}^-(\lambda^-)\varepsilon^kx^i(-y)^j \\[2.5ex]
x-\vartheta_1y+\displaystyle\sum\limits_{k-1,i+j=2}^{n}
\widetilde{B_{kij}}^-(\lambda^-)\varepsilon^kx^i(-y)^j
\end{array}\right), \ \ &\text{if} \ \ y>0,\\[1.0ex]
&\left(\begin{array}{c}
-\vartheta_2x-y+\displaystyle\sum\limits_{k=1,i+j=2}^{n}
\widetilde{A_{kij}}^+(\lambda^+)\varepsilon^kx^i(-y)^j \\[2.5ex]
x-\vartheta_1y+\displaystyle\sum\limits_{k=1,i+j=2}^{n}
\widetilde{B_{kij}}^+(\lambda^+)\varepsilon^kx^iy^j
\end{array}\right), \ \ &\text{if} \ \ y\leq0.
\end{aligned}\right.
\end{equation}
Then, the displacement function can be written as
\begin{equation}\label{Eqn2-9}
\begin{aligned}
d(\xi)=\Pi^+(\xi)-(\Pi^-)^{-1}(\xi)=\Pi^+(\xi)-\Pi_-^+(\xi)
=\sum\limits_{j\geq1}V_j(\lambda^{\pm},\varepsilon,\pi)\xi^j,
\end{aligned}
\end{equation}where $\Pi_-^+(\xi)$ is the first half-return map of \eqref{Eqn2-7-0}, and
\begin{equation}\label{Eqn2-10}
V_j(\lambda^{\pm},\varepsilon,\pi)=\sum_{k=1}^\infty
V_{jk}\, \varepsilon^k , \quad j=1,2,\ldots,
\end{equation}
in which $V_{jk}$ denotes the $j$th $\varepsilon^k$-order Lyapunov constant,
see \cite{Chen4} for more details about the computation of the generalized
Lyapunov constants.
Hence, we can derive the center conditions of system \eqref{Eqn2-4}
by vanishing the $\varepsilon$ terms in these generalized Lyapunov constants,
and then by Lemma \ref{L0} we derive the algebraic conditions
characterizing the nilpotent center of system \eqref{Eqn2-1}.
Thus, we prove that these conditions are necessary
for $(0,0)$ of \eqref{Eqn2-1} be a nilpotent center. In general,
these nilpotent center conditions can be satisfied by appropriately choosing
the perturbation coefficients $P_{kij}^{\pm}$ and $Q_{kij}^{\pm}$.

Further, we generalize our method to study the bifurcation of limit cycles
from the nilpotent center of the switching system \eqref{Eqn2-1}.
Actually, by B-T bifurcation theory, we may add the linear perturbation
term $-\varepsilon^2 x$ to the switching system \eqref{Eqn2-1} such that the system has
a linear-type center at the origin.
Then we perturb such system by adding higher
$\varepsilon$-order terms,}
and compute the generalized
Lyapunov constants of the perturbed system.
Finally, by using the higher $\varepsilon$-order Lyapunov constants
\cite{TY2018,YHL2018}, we find the bifurcation of limit cycles
from the center as many as possible.
More precisely, with the result in \cite{TIAN}, we derive
the following lemma giving the sufficient conditions
for the existence of small-amplitude limit cycles
around the origin of \eqref{Eqn2-1}.

\begin{lemma}[\rm [Lemma 4, \cite{TIAN}]\label{Lem2.3}
If there exists a  critical point
$\lambda_*=(a_{1c},a_{2c},\cdots,a_{nc})$ such that a
set of $\varepsilon^k$-order Lyapunov constants satisfies
$V_{1k}(\lambda_*)=V_{2k}(\lambda_*)=\cdots=V_{nk}(\lambda_*)=0$,
$V_{n+1,k}(\lambda_*)\neq0$ and
\begin{equation}\label{Eqn2-11}
\begin{array}{l}
{\rm{det}}\bigg[\displaystyle\frac{\partial(V_{1k},V_{2k},\cdots,
V_{nk})}{\partial(a_{1c},a_{2c},\cdots,a_{nc})} (\lambda_*)\bigg]\neq0,
\end{array}
\end{equation}
then small appropriate perturbations about $\lambda=\lambda_*$
lead to that the switching system \eqref{Eqn2-7} has exactly $n$ limit cycles bifurcating
from $(0,0)$.
\end{lemma}

\section{The proof of Theorem \ref{Th1}}

From the results given in the precious section,
we know that the nilpotent center of \eqref{Eqn1-5} at the origin
can be combined by a monodromic singular point or a cusp,
i.e. the parameters satisfy the condition $\Xi_1$ or $\Xi_2$.
Firstly, assume that the origin in the
first smooth system of \eqref{Eqn1-5} is a cusp, then we have $b_2^+\neq0$.
If $b_2^-\neq0$, then the origin in the second system
of \eqref{Eqn1-5} is also a cusp with multiplicity two.
If $b_2^-=0$, we have the following statement: The origin of the second
system of \eqref{Eqn1-5} is a monodromic singular point with
multiplicity three if and only if one of the following conditions holds:
\begin{equation}\label{Eqn3-1}
\begin{aligned}
&a_2^-=b_2^-=0, \quad b_3^-> \dfrac{3}{4}\, (a_3^-)^2;\\
&b_2^-=0, \quad a_2^-\neq0, \quad b_3^-> \dfrac{1}{3}\, (a_2^-)^2.
\end{aligned}
\end{equation}
We only consider $b_2^+>0$ and $b_2^-\leq 0$ because the origin of
\eqref{Eqn1-5} cannot be monodromic when either $b_2^+< 0$ or $b_2^-> 0$,
which are illustrated in the
two phase portraits of system \eqref{Eqn1-5}, as
depicted in Figure \ref{Fig1}.
It is shown that when $b_2^+<0$ the first system of
\eqref{Eqn1-5} has two seperatrices connecting the origin and
two singular points in the right half Poincar\'e disc, see
Figure \ref{Fig1a}. Similarly when $b_2^->0$ the second system
of \eqref{Eqn1-5} has two seperatrices connecting the origin
and the other two ones in the left half Poincar\'e disc,
see Figure \ref{Fig1b}. Thus, the origin cannot be a center
when either $b_2^+< 0$ or $b_2^-> 0$.

\begin{example}
Two different phase portraits of the cubic switching Li\'enard
system \eqref{Eqn1-5}, as depicted in the Poincar\'e disc,
show that the origin is a cusp for $b_2^+=-1$ (see Figure \ref{Fig1}(a))
or $b_2^-=3$ (see Figure \ref{Fig1}(b)).
\end{example}

\begin{figure}[!t]
\vspace{0.00in}
\centering
\hspace*{-0.0in}
\subfigure[]{
\label{Fig1a}
\includegraphics[width=0.28\textwidth,height=0.20\textheight]{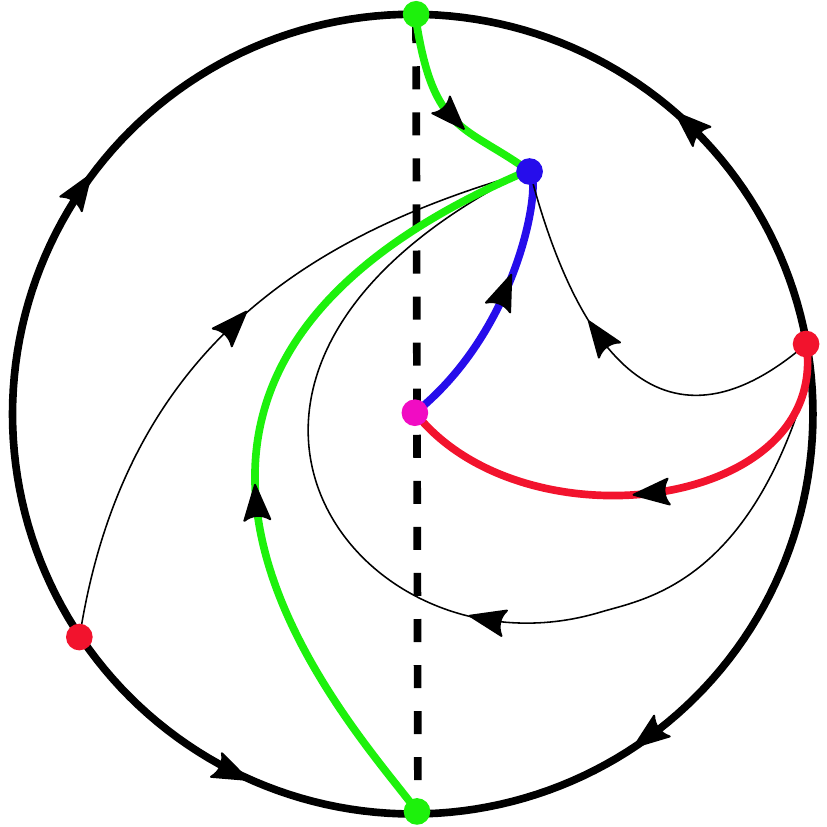}}
\hspace*{0.80in}
\subfigure[]{
\label{Fig1b}
\includegraphics[width=0.28\textwidth,height=0.20\textheight]{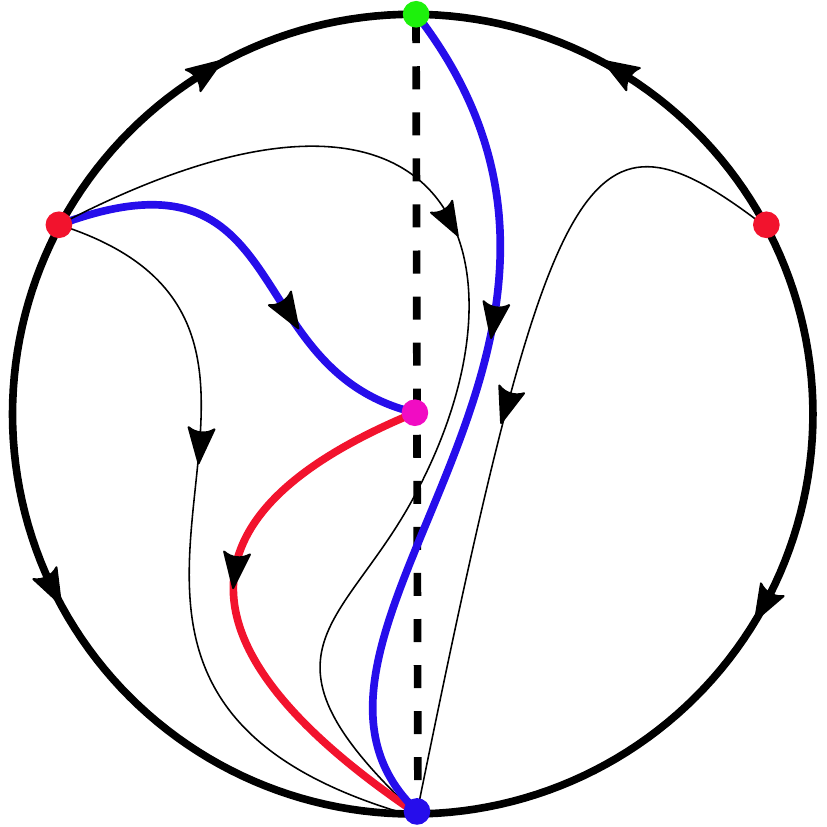}}
\hspace*{0.15in}
\caption{Two phase portraits of the cubic switching Li\'enard
system \eqref{Eqn1-5} showing that the origin is a cusp for
(a) $a_2^{\pm}=b_{3}^{\pm}=1$, $a_3^{\pm}=2$, $b_2^+=-1$ and $b_2^-=-2$;
and (b) $a_2^{\pm}=a_3^{\pm}=b_{3}^+=1$, $b_3^{-}=-1$, $b_2^+=2$
and $b_2^-=3$.}
\label{Fig1}
\end{figure}

Next, we discuss how to apply the higher-order Poincar\'e-Lyapunov method to
derive the nilpotent center conditions for the origin of system \eqref{Eqn1-5}.
With perturbations added
to system \eqref{Eqn1-5}, we obtain the following perturbed systems,
\begin{equation}\label{Eqn3-3}
\left(\begin{array}{cc}
\dot{x}\\
\dot{y}
\end{array}\right)
=\left\{\begin{aligned}&\left(\begin{array}{c}
y-(a_2^+x^2+a_3^+x^3)+\displaystyle\sum_{k=1}\varepsilon^k P_{k}^+(x) \\[3.0ex]
-\varepsilon^2 x-(b_2^+x^2+b_3^+x^3)+\displaystyle\sum_{k=1}\varepsilon^k
Q_{k}^+(x)
\end{array}\right), &\text{if} \ \ x\geq0,\\[1.0ex]
&\left(\begin{array}{c}
y-(a_2^-x^2+a_3^-x^3)+\displaystyle\sum_{k=1}\varepsilon^k P_{k}^-(x) \\[3.0ex]
-\varepsilon^2 x-(b_2^-x^2+b_3^-x^3)+\displaystyle\sum_{k=1}\varepsilon^k
Q_{k}^-(x)
\end{array}\right), &\text{if} \ \ x<0,
\end{aligned}
\right.
\end{equation}
where
\begin{equation*}
P_{k}^{\pm}(x)=p_{k2}^{\pm}\, x^2+p_{k3}^{\pm}\, x^3,\quad
Q_{k}^{\pm}(x)=q_{k2}^{\pm}\, x^2+q_{k3}^{\pm}\, x^3,
\end{equation*}
in which $p_{ki}^{\pm}$, $q_{ki}^{\pm}$ are real parameters.
It is very difficult for computing the generalized Lyapuov constants
of \eqref{Eqn3-3} with large number of parameters.
So we only consider the $\varepsilon^2$-order perturbations.
Further, introducing the transformation
$(x,y,t)\rightarrow (\varepsilon^3 x, \varepsilon^2 y, \frac{t}{\varepsilon})$
into system \eqref{Eqn3-3}, we obtain
\begin{equation}\label{Eqn3-4}
\left(\begin{array}{cc}
\dot{x}\\
\dot{y}
\end{array}\right)
=\left\{\begin{aligned}&\left(\begin{array}{c}
-y-(b_2^++q_{22}^+ \varepsilon^2)y^2-(b_3^+\varepsilon^2
+ q_{23}^+\varepsilon^4)y^3 \\[1.0ex]
x-(a_2^+\varepsilon+ p_{22}^+\varepsilon^3) y^2
-(a_3^+\varepsilon^3+p_{23}^+\varepsilon^5)y^3
\end{array}\right), &\text{if} \ \ y\geq0,\\[1.0ex]
&\left(\begin{array}{c}
-y-(b_2^-+q_{22}^- \varepsilon^2)y^2-(b_3^-\varepsilon^2
+q_{23}^-\varepsilon^4)y^3 \\[1.0ex]
x-(a_2^-\varepsilon+p_{22}^-\varepsilon^3) y^2-(a_3^-\varepsilon^3
+p_{23}^-\varepsilon^5)y^3
\end{array}\right), &\text{if} \ \ y<0.
\end{aligned}
\right.
\end{equation}
Then, we consider two cases: $b_2^-=0$ and $b_2^-<0$ to obtain
the nilpotent center conditions by considering
each $i$th $\varepsilon^k$-order Lyapunov constant.

\subsection{Case 1: $b_2^-=0$}

We use the higher-order
Poincar\'{e}-Lyapunov method to compute the generalized
Lyapunov constants associated with the origin of \eqref{Eqn3-4}.
The first two generalized Lyapunov constants are $V_1(\varepsilon)=0$ and
\begin{equation}\label{V2}
V_2(\varepsilon)= \frac{4}{3}\big[(a_2^+-a_2^-)
+(p_{22}^+-p_{22}^-)\varepsilon^2 \big]\, \varepsilon.
\end{equation}
By \eqref{Eqn3-1} we consider the two subcases:
(i) $a_2^-=0$, $b_3^-> \frac{3}{4} (a_3^-)^2$,
and (ii) $a_2^-\neq0$, $b_3^-> \frac{1}{3} (a_2^-)^2$.

\begin{enumerate}
\item[{(i)}]
$a_2^-=0$, $b_3^-> \frac{3}{4} (a_3^-)^2$.
Setting the $\varepsilon$-order and $\varepsilon^3$-order
Lyapunov constants in $V_2(\varepsilon)$ zero yields
the necessary center conditions $a_2^+=0$ and $p_{22}^-=p_{22}^+$.
Then, the $3$rd generalized Lyapunov constant is obtained as
\begin{equation*}
V_3(\varepsilon)= -\, \frac{\pi}{8}
\big[3a_3^+ + 3a_3^- - 2b_2^+p_{22}^+ + (3p_{23}^++3p_{23}^-
-2p_{22}^+q_{22}^+
-2p_{22}^+q_{22}^-)\varepsilon^2 \big]\, \varepsilon^3.
\end{equation*}
Letting the $\varepsilon^{3}$-order and $\varepsilon^5$-order terms
in $V_3(\varepsilon)$ equal zero we obtain the following conditions,
\begin{equation*}
\begin{aligned}
p_{22}^+=&\ \frac{3(a_3^++a_3^-)}{2b_2^+},\\
p_{23}^-=&\ \frac{1}{b_2^+}\big(\! -b_2^+ p_{23}^+
+a_3^+q_{22}^+ +a_3^-q_{22}^++a_3^+q_{22}^-+a_3^-q_{22}^- \big).
\end{aligned}
\end{equation*}
Then, we have the 4th generalized Lyapunov constant, given by
\begin{equation*}
V_4(\varepsilon)=\frac{4}{15 b_2^+}\,
\big[4 a_3^- (b_2^+)^2-M_{1}\varepsilon^2+M_{2}\varepsilon^4\big]\,
\varepsilon^3,
\end{equation*}where
\begin{equation*}
\begin{aligned}
M_{1}=&\ 3a_3^+b_3^++3a_3^-b_3^+-3a_3^+b_3^--3a_3^-b_3^--4(b_2^+)^2p_{23}^-
-4a_3^-b_2^+q_{22}^++4a_3^+b_2^+q_{22}^-,\\[0.5ex]
M_{2}=&\ 4b_2^+p_{23}^-q_{22}^++4b_2^+p_{23}^-q_{22}^--4a_3^+q_{22}^+q_{22}^-
-4a_3^-q_{22}^+q_{22}^--4a_3^+(q_{22}^-)^2-4a_3^-(q_{22}^-)^2\\
&-3a_3^+q_{23}^+-3a_3^-q_{23}^++3a_3^+q_{23}^-+3a_3^-q_{23}^-.\\
\end{aligned}
\end{equation*}
Letting each $\varepsilon$ term in $V_4(\varepsilon)$ be zero,
we obtain the conditions:
\begin{equation*}
\begin{aligned}
a_3^-=&\ 0,\\
p_{23}^-=&\ \frac{3a_3^+b_3^+-3a_3^+b_3^-+4a_3^+b_2^+q_{22}^-}{4(b_2^+)^2},\\
q_{23}^-=&\ \frac{1}{b_2^+}\big(\! -b_3^+q_{22}^++b_3^-q_{22}^+
-b_3^+q_{22}^-+b_3^-q_{22}^-+b_2^+q_{23}^+\big).
\end{aligned}
\end{equation*}
Then, we have
\begin{equation*}
V_5(\varepsilon)=-\frac{5\pi}{256 (b_2^+)^3}a_3^+\,(b_2^+
+q_{22}^+\varepsilon^2+q_{22}^-\varepsilon^2)
\big[2(b_2^+)^2(b_3^++7b_3^-)+M_{3}\varepsilon^2\big]\, \varepsilon^5,
\end{equation*}
where
\begin{equation*}
\begin{aligned}
M_{3}=&\ 9(b_3^+)^2-18b_3^+b_3^-+9(b_3^-)^2-14b_2^+b_3^+q_{22}^+
+14b_2^+b_3^-q_{22}^+-2b_2^+b_3^+q_{22}^-\\
&+2b_2^+b_3^-q_{22}^-+16(b_2^+)^2q_{23}^+.
\end{aligned}
\end{equation*}
$a_3^+=0$ (leading to $V_5(\varepsilon)=0$) gives the condition I.
Otherwise, if $a_3^+ \ne 0$,
setting the other $\varepsilon$ terms in $V_5(\varepsilon)$ zero we obtain
\begin{equation*}
\begin{aligned}
b_3^+=-\,7b_3^-, \quad
q_{23}^+= -\,\frac{1}{(b_2^+)^2} \big[
36(b_3^-)^2 + 7b_2^+b_3^-q_{22}^+ + b_2^+b_3^-q_{22}^- \big].
 \end{aligned}
\end{equation*}
Further, under the above conditions, we have the 6th generalized Lyapunov
constant given by
\begin{equation*}
\begin{aligned}
V_6(\varepsilon)=&-\frac{128}{315 (b_2^+)^2}a_3^+b_3^-\,
(b_2^++q_{22}^+\varepsilon^2+q_{22}^-\varepsilon^2)^3\varepsilon^5\neq0
\end{aligned}
\end{equation*}when $a_3^+\neq0$.
\end{enumerate}

\begin{enumerate}
\item[{(ii)}]
$a_2^-\neq0$, $b_3^-> \frac{1}{3} (a_2^-)^2$.
Setting the $\varepsilon$-order and $\varepsilon^3$-order terms
in $V_2(\varepsilon)$ zero yields the necessary center conditions
$a_2^-=a_2^+\neq 0$ and $p_{22}^-=p_{22}^+$.
Then, the $3$rd generalized Lyapunov constant is given by
\begin{equation*}
\begin{aligned}
V_3(\varepsilon)=&\ \frac{\pi}{8}
\big[2a_2^-b_2^+-(3a_3^++3a_3^--2b_2^+p_{22}^+-2a_2^+q_{22}^+
-2a_2^+q_{22}^-)\varepsilon^2-(3p_{23}^++3p_{23}^-\\
&-2p_{22}^+q_{22}^+-2p_{22}^+q_{22}^-)\varepsilon^4 \big]\varepsilon\neq0.
\end{aligned}
\end{equation*}
\end{enumerate}

\subsection{Case 2: $b_2^-<0$}

Similarly, setting each $\varepsilon$ term in
$V_2(\varepsilon)$ zero yields the necessary center conditions
$a_2^-=a_2^+$ and $p_{22}^-=p_{22}^+$.
Then, we obtain the $3$rd generalized Lyapunov constant,
\begin{equation*}
\begin{aligned}
V_3(\varepsilon)=&\ \frac{\pi}{8}
\big[2a_2^+(b_2^++b_2^-)-(3a_3^++3a_3^--2b_2^+p_{22}^+-2b_2^-p_{22}^+-2a_2^+q_{22}^+-2a_2^+q_{22}^-)\varepsilon^2\\
&-(3p_{23}^++3p_{23}^--2p_{22}^+q_{22}^+-2p_{22}^+q_{22}^-)\varepsilon^4 \big]\varepsilon.
\end{aligned}
\end{equation*}
Considering the $\varepsilon$-order term in $V_3(\varepsilon)$,
we have two subcases: (i) $a_2^+=0$, and (ii) $b_2^++b_2^-=0$.

\begin{enumerate}
\item[{(i)}]
Assume that $a_2^+=0$ and $b_2^++b_2^-\neq0$. We let the $\varepsilon^{3}$-order and $\varepsilon^5$-order terms
in $V_3(\varepsilon)$ equal zero to obtain the conditions,
\begin{equation*}
\begin{aligned}
p_{22}^+=&\ \frac{3(a_3^++a_3^-)}{2(b_2^++b_2^-)},\\
p_{23}^-=&\ \frac{1}{b_2^++b_2^-}\big(b_2^+p_{23}^+ + b_2^-p_{23}^+
- a_3^+q_{22}^+ - a_3^-q_{22}^+ - a_3^+q_{22}^- - a_3^-q_{22}^- \big).
\end{aligned}
\end{equation*}
Then, the 4th generalized Lyapunov constant becomes
\begin{equation*}
\begin{aligned}
V_4(\varepsilon)=&\ \frac{-\,4}{15 (b_2^++b_2^-)}\,
\big[4(b_2^++b_2^-)(a_3^-b_2^+-a_3^+b_2^-)
-M_4\varepsilon^2-M_5\varepsilon^4\big]\varepsilon^3,
\end{aligned}
\end{equation*}where
\begin{equation*}
\begin{aligned}
M_4=&\ 3a_3^+b_3^++3a_3^-b_3^+-3a_3^+b_3^--3a_3^-b_3^-+4(b_2^+)^2p_{23}^+
+8b_2^+b_2^-p_{23}^++4(b_2^-)^2p_{23}^+ \\
&-4a_3^+b_2^+q_{22}^+-8a_3^-b_2^+q_{22}^+-4a_3^-b_2^-q_{22}^+
-4a_3^-b_2^+q_{22}^-+4a_3^+b_2^-q_{22}^-,\\[0.5ex]
M_5=&\ 4b_2^+p_{23}^+q_{22}^++4b_2^-p_{23}^+q_{22}^+-4a_3^+(q_{22}^+)^2
-4a_3^-(q_{22}^+)^2+4b_2^+p_{23}^+q_{22}^-+4b_2^-p_{23}^+q_{22}^- \\
&-4a_3^+q_{22}^+q_{22}^--4a_3^-q_{22}^+q_{22}^-+3a_3^+q_{23}^++3a_3^-q_{23}^+-3a_3^+q_{23}^--3a_3^-q_{23}^-.
\end{aligned}
\end{equation*}
Thus, we have
\begin{equation*}
a_3^-=\frac{a_3^+b_2^-}{b_2^+},\quad
p_{23}^+=\frac{-3a_3^+b_2^++3a_3^+b_3^-+4a_3^+b_2^+q_{22}^+
+4a_3^+b_2^-q_{22}^+}{4b_2^+(b_2^++b_2^-)}
\end{equation*}
by setting the $\varepsilon^{3}$-order and $\varepsilon^{5}$-order terms
in $V_4(\varepsilon)$ zero. Then, $M_5$ is simplified as
\begin{equation*}
M_5=-\frac{3}{b_2^+}a_3^+(b_3^+q_{22}^+-b_3^-q_{22}^++b_3^+q_{22}^-
-b_3^-q_{22}^--b_2^+q_{23}^+-b_2^-q_{23}^++b_2^+q_{23}^-+b_2^-q_{23}^-).
\end{equation*}

(i.1) If $a_3^+=0$, we obtain $M_5=0$,
leading to a condition included in the condition II.

(i.2) If $a_3^+ \ne 0$, we have
$$
q_{23}^+=\frac{1}{b_2^++b_2^-}(b_3^+q_{22}^+-b_3^-q_{22}^++b_3^+q_{22}^-
-b_3^-q_{22}^-+b_2^+q_{23}^-+b_2^-q_{23}^-)
$$
by setting $M_5=0$. Further, we obtain the $5$th generalized Lyapunov constant,
\begin{equation*}
\begin{aligned}
V_5(\varepsilon)=&-\frac{5a_3^+\pi}{256b_2^+(b_2^++b_2^-)^2} \,
\big[b_2^++b_2^-+(q_{22}^++q_{22}^-)\varepsilon^2\big]
\big[2(b_2^++b_2^-)M_6+M_7\varepsilon^2\big]\varepsilon^5
\end{aligned}
\end{equation*}
where
\begin{equation*}
\begin{aligned}
M_6=&\ b_3^+(b_2^++7b_2^-)+b_3^-(7b_2^++b_2^-),\\[0.5ex]
M_7=&\ 9(b_3^+)^2-18b_3^+b_3^-+9(b_3^-)^2+2b_2^+b_3^+q_{22}^+
+2b_2^-b_3^+q_{22}^+-2b_2^+b_3^-q_{22}^+-2b_2^-b_3^-q_{22}^+\\
&+14b_2^+b_3^+q_{22}^-+14b_2^-b_3^+q_{22}^--14b_2^+b_3^-q_{22}^-
-14b_2^-b_3^-q_{22}^-+16(b_2^+)^2q_{23}^-\\
&+32b_2^+b_2^-q_{23}^-+16(b_2^-)^2q_{23}^-.
\end{aligned}
\end{equation*}
Since $a_3^+(b_2^++b_2^-)\neq0$,
the only possibility for $V_5(\varepsilon)=0$
is $M_6=M_7=0$. Considering $M_6=0$, if $b_2^++7b_2^-\neq0$,
we have $b_3^+=\frac{-(7b_2^++b_2^-)b_3^-}{b_2^++7b_2^-}$;
and if $b_2^+=-7b_2^-$, we obtain $M_6=-48b_2^-b_3^-=0$ which
yields $b_3^-=0$ due to $b_2^-<0$.
Thus, we consider the following two subcases.

(i.2.1) If $b_3^+=\frac{-(7b_2^++b_2^-)b_3^-}{b_2^++7b_2^-}$, setting
$M_7=0$ yields
\begin{equation*}
\begin{aligned}
q_{23}^-=\frac{1}{(b_2^++7b_2^-)^2}[-36(b_3^-)^2+b_2^+b_3^-q_{22}^++7b_2^-b_3^-q_{22}^++7b_2^+b_3^-q_{22}^-+49b_2^-b_3^-q_{22}^-].
\end{aligned}
\end{equation*}
Then, we obtain the $6$th generalized Lyapunov constant,
\begin{equation*}
\begin{aligned}
V_6(\varepsilon)=&-\frac{128}{315 b_2^+(b_2^++7b_2^-)}a_3^+b_3^-
\big[b_2^++b_2^-+(q_{22}^++q_{22}^-)\varepsilon^2\big]^3\varepsilon^5.
\end{aligned}
\end{equation*}
If $b_3^-=0$, we have $V_6(\varepsilon)=0$,
and this necessary condition is included in the condition III.
Otherwise, we have  $V_6(\varepsilon)\neq0$ when  $b_2^++b_2^-\neq0$.

(i.2.2) If $b_3^-=b_2^++7b_2^-=0$  we have
\begin{equation*}
\begin{aligned}
q_{23}^-=&\frac{1}{192(b_2^-)^2}[-3(b_3^+)^2+4b_2^-b_3^+q_{22}^-+28b_2^-b_3^+q_{22}^-]
\end{aligned}
\end{equation*}from $M_7=0$. Then the $6$th generalized Lyapunov constant has the following form
\begin{equation*}
\begin{aligned}
V_6(\varepsilon)=&-\frac{8}{6615(b_2^-)^2}a_3^+b_3^+
\big[6b_2^-+(q_{22}^++q_{22}^-)\varepsilon^2\big]^3\varepsilon^5\neq0.
\end{aligned}
\end{equation*}Assume that $b_3^+=0$, we have $V_6(\varepsilon)=0$.
Combining the necessary conditions from (i.2.1) and (i.2.2),
we obtain the condition III. Otherwise we have
$V_6(\varepsilon)\neq0$ when $a_3\neq0$.
\end{enumerate}

\begin{enumerate}
\item[{(ii)}]
Assume that $b_2^-+b_2^+=0$. Consider two subcases:
(ii.1) $a_2^+=0$ and (ii.2) $a_2^+\neq0$.

(ii.1) If $a_2^+=0$,
setting each $\varepsilon$ term in $V_3(\varepsilon)$ zero we have
\begin{equation*}
\begin{aligned}
a_3^-=&-a_3^+,\\
p_{23}^-=&\ \frac{1}{3}(-3p_{23}^++2p_{22}^+q_{22}^++2p_{22}^+q_{22}^-).
\end{aligned}
\end{equation*}
Then, the 4th generalized Lyapunov constant is given by
\begin{equation*}
\begin{aligned}
V_4(\varepsilon)=&-\frac{8}{45}\,
\big(M_8+M_9\varepsilon^2\big)\, \varepsilon^5,
\end{aligned}
\end{equation*}where
\begin{equation*}
\begin{aligned}
M_8=&-3b_3^+p_{22}^++3b_3^-p_{22}^+-6a_3^+q_{22}^++4b_2^+p_{22}^+q_{22}^+
-6a_3^+q_{22}^-+4b_2^+p_{22}^+q_{22}^-,\\[0.5ex]
M_9=&\ 6p_{23}^+q_{22}^+-4p_{22}^+(q_{22}^+)^2+6p_{23}^+q_{22}^-
-4p_{22}^+q_{22}^+q_{22}^-+3p_{22}^+q_{23}^+-3 p_{22}^+q_{23}^-.\\
\end{aligned}
\end{equation*}
(ii.1.1)
If $q_{22}^-=-q_{22}^+$, it follows from $M_8=M_9=0$ that
either $b_3^-=b_3^+$, $q_{23}^-=q_{23}^+$, or $p_{22}=0$.
The first choice leads to a necessity condition included in the
condition IV. For the second choice, we have the 5th generalized
Lyapunov constant,
\begin{equation*}
\begin{aligned}
V_5(\varepsilon)=&\ \frac{15}{64}\,
\big(a_3^++p_{23}^+\varepsilon^2\big)\big(b_3^+-b_3^-+q_{23}^+
\varepsilon^2-q_{23}^-\varepsilon^2\big)\, \varepsilon^5,
\end{aligned}
\end{equation*}
which yields $a_3^+=p_{23}^+=0$ by setting $V_5(\varepsilon)=0$,
leading to a condition included in the condition II.
Combining the subcases (i.1) and (ii.1.1) we have the condition II.

(ii.1.2) If $q_{22}^-\neq -q_{22}^+$, we have
\begin{equation*}
\begin{aligned}
a_3^+=&\ \frac{1}{6(q_{22}^++q_{22}^-)}(-3 b_3^+ p_{22}^+ + 3 b_3^- p_{22}^+
+ 4 b_2^+ p_{22}^+ q_{22}^+ + 4 b_2^+ p_{22}^+ q_{22}^-),\\[1.0ex]
p_{23}^+=&\ \frac{1}{6(q_{22}^++q_{22}^-)}\big[4 p_{22}^+ (q_{22}^+)^2
+ 4 p_{22}^+ q_{22}^+ q_{22}^- - 3 p_{22}^+ q_{23}^+ + 3 p_{22}^+ q_{23}^-
\big]
\end{aligned}
\end{equation*}
by setting $M_8=M_9=0$. Then, we obtain the 5th generalized Lyapunov constant,
\begin{equation*}
\begin{aligned}
V_5(\varepsilon)=&-\frac{5\pi}{384(q_{22}^++q_{22}^-)}p_{22}^+
\big(M_{10}+2M_{11}\varepsilon^2+M_{12}\varepsilon^2\big)\, \varepsilon^5,
\end{aligned}
\end{equation*}where
\begin{equation*}
\begin{aligned}
M_{10}=&\ 3 (b_3^+ - b_3^-) (3 b_3^+ - 3 b_3^- - 4 b_2^+ q_{22}^+
- 4 b_2^+ q_{22}^-),\\[0.5ex]
M_{11}=&\ b_3^+ (q_{22}^+)^2 + 7 b_3^- (q_{22}^+)^2
+ 8 b_3^+ q_{22}^+ q_{22}^- + 8 b_3^- q_{22}^+ q_{22}^-
+ 7 b_3^+ (q_{22}^-)^2 + b_3^- (q_{22}^-)^2+ 9 b_3^+ q_{23}^+ \\
&- 9 b_3^- q_{23}^+ - 6 b_2^+ q_{22}^+ q_{23}^+ - 6 b_2^+ q_{22}^- q_{23}^+
- 9 b_3^+ q_{23}^- + 9 b_3^- q_{23}^- + 6 b_2^+ q_{22}^+ q_{23}^-
+ 6 b_2^+ q_{22}^- q_{23}^-,\\[0.5ex]
M_{12}=&\ 2 (q_{22}^+)^2 q_{23}^+ + 16 q_{22}^+ q_{22}^- q_{23}^+
+ 14 (q_{22}^-)^2 q_{23}^+ + 9 (q_{23}^+)^2 +
 14 (q_{22}^+)^2 q_{23}^- + 16 q_{22}^+ q_{22}^- q_{23}^- \\
 &+ 2 (q_{22}^-)^2 q_{23}^- - 18 q_{23}^+ q_{23}^- + 9 (q_{23}^-)^2.
\end{aligned}
\end{equation*}
If $p_{22}^+=0$ or $b_2^+=\frac{3(b_3^--b_3^+)}{4(q_{22}^++q_{22}^-)}$,
we have $a_3^+=0$, which is included in the condition II.
Otherwise, we have $b_3^-=b_3^+$ from $M_{10}$,
which leads to a necessity condition included in the condition IV.

(ii.2) If $a_2^+ \ne 0$, letting the $\varepsilon^{3}$-order
and $\varepsilon^5$-order terms
in $V_3(\varepsilon)$ equal zero we obtain the conditions:
\begin{equation*}
\begin{aligned}
q_{22}^-=&\ \frac{3a_3^++ 3a_3^--2a_2^+q_{22}^+}{2a_2^+},\\
p_{23}^-=&-\frac{a_3^+p_{22}^++a_3^-p_{22}^+-a_2^+p_{23}^+}{a_2^+}.
\end{aligned}
\end{equation*}
Then, we have the 4th generalized Lyapunov constant, given by
\begin{equation*}
\begin{aligned}
V_4(\varepsilon)=&\ \frac{8}{15 a_2^+}\,
\big(a_2^+M_{13}+M_{14}\varepsilon^2+M_{15}\varepsilon^4\big)\, \varepsilon^3,
\end{aligned}
\end{equation*}
where
\begin{equation*}
\begin{aligned}
M_{13}=&-2a_3^+b_2^+-2a_3^-b_2^++a_2^+b_3^+-a_2^+b_3^-,\\[0.5ex]
M_{14}=&\ 3(a_3^+)^2+3a_3^+a_3^--2a_3^+b_2^+p_{22}^+-2a_3^-b_2^+p_{22}^+
+a_2^+b_3^+ p_{22}^+-a_2^+b_3^-p_{22}^+\\
 &-2a_2^+a_3^+q_{22}^+-2a_2^+a_3^-q_{22}^++(a_2^+)^2q_{23}^+
-(a_2^+)^2q_{23}^-,\\[0.5ex]
M_{15}=&\ 3a_3^+p_{23}^++3a_3^-p_{23}^+-2a_3^+p_{22}^+q_{22}^+
-2a_3^-p_{22}^+q_{22}^++a_2^+p_{22}^+q_{23}^+-a_2^+p_{22}^+q_{23}^-.\\
\end{aligned}
\end{equation*}
Setting $M_{13}=M_{14}=0$, we obtain the conditions:
\begin{equation*}
\begin{aligned}
b_3^-=&\ \frac{1}{a_2^+}(-2a_3^+b_2^+-2a_3^-b_2^++a_2^+b_3^+),\\
q_{23}^-=&\ \frac{1}{(a_2^+)^2} \big[3(a_3^+)^2+3a_3^+a_3^-
-2a_2^+a_3^+q_{22}^+-2a_2^+a_3^-q_{22}^++(a_2^+)^2q_{23}^+ \big].\\
\end{aligned}
\end{equation*}
Then, two subcases follow from $M_{15}=0$:
(ii.2.1) $a_3^++a_3^-=0$ and (ii.2.2) $p_{23}^+=\frac{a_3^+p_{22}^+}{a_2^+}$.

(ii.2.1) If $a_3^++a_3^-=0$, which is combined with the necessary conditions from (ii.1.1)
to yield the condition IV.

(ii.2.2) If $p_{23}^+=\frac{a_3^+p_{22}^+}{a_2^+}$, then
we have $V_4(\varepsilon)=0$ and obtain the $5$th generalized Lyapunov
constant,
\begin{equation*}
\begin{aligned}
V_5(\varepsilon)=&-\frac{5\pi}{128 (a_2^+)^3} \,(a_3^+ + a_3^-)(a_2^++p_{22}^+\varepsilon^2)
\big[2 a_2^+ M_{16}+M_{17}\varepsilon^2\big]\, \varepsilon^5,
\end{aligned}
\end{equation*}
where
\begin{equation*}
\begin{aligned}
M_{16}=&-7a_3^+ b_2^+-a_3^-b_2^++4a_2^+b_3^+,\\[0.5ex]
M_{17}=&\ 21(a_3^+)^2+3a_3^+a_3^--14a_2^+a_3^+q_{22}-2a_2^+a_3^-q_{22}^+
+8(a_2^+)^2q_{23}^+.
\end{aligned}
\end{equation*}
Setting $M_{16}=M_{17}=0$ results in
\begin{equation*}
\begin{aligned}
b_3^+=&\ \frac{7a_3^+b_2^++a_3^-b_2^+}{4a_2^+},\\
q_{23}^+=&- \frac{1}{q_{23}^+} \big[21(a_3^+)^2+3a_3^+a_3^-
-14a_2^+a_3^+q_{22}^+ - 2a_2^+a_3^-q_{22}^+ \big],
\end{aligned}
\end{equation*}
under which the $6$th generalized Lyapunov constant becomes
\begin{equation*}
\begin{aligned}
V_6(\varepsilon)=&\ \frac{8}{105 (a_2^+)^4}(a_3^+ + a_3^-)^3
(a_2^++p_{22}^+\varepsilon^2)(2a_2^+b_2^+-3a_3^+\varepsilon^2
+2 a_2^+q_{22}^+\varepsilon^2)\, \varepsilon^7\neq0
\end{aligned}
\end{equation*}
with $a_3^++a_3^-\neq0$.
\end{enumerate}

The above results give the necessary nilpotent center conditions
for system \eqref{Eqn1-5} at the origin when the origin of
the first system of \eqref{Eqn1-5} is a cusp.
Now we consider the center conditions associated with the origin of the
first system of \eqref{Eqn1-5} being a monodromic singular point
with multiplicity three if and only if one of the following holds:
\begin{equation}\label{Eqn3-5}
\begin{aligned}
&a_2^+=b_2^+=0, \quad b_3^+> \dfrac{3}{4}\, (a_3^+)^2;\\
&b_2^+=0, \quad a_2^+\neq0, \quad b_3^+> \dfrac{1}{3}\, (a_2^+)^2.
\end{aligned}
\end{equation}
If $b_2^-<0$, the origin in the second system of \eqref{Eqn1-5} is a cusp,
which is a similar case to that discussed above, leading to
the condition V in Theorem \ref{Th1}.
Hence, we assume that one of the conditions in \eqref{Eqn3-1} holds.
Then, the origin in the second system of \eqref{Eqn1-5} is also
a monodromic singular point.

Setting the second generalized Lyapunov constant
\eqref{V2} zero,
we find two cases: either $a_2^+=a_2^-=0$ or $a_2^+a_2^-\neq0$ when the
conditions in \eqref{Eqn3-1} and \eqref{Eqn3-5} are satisfied.
The detailed analysis for the two cases are given below.

\subsection{Case 3: $a_2^+=a_2^-=0$}

The second generalized Lyapunov constant becomes
\begin{equation*}
\begin{aligned}
V_2(\varepsilon)=&\ \frac{4}{3}\,(p_{22}^+-p_{22}^-)\, \varepsilon^3.
\end{aligned}
\end{equation*}
Setting $V_2(\varepsilon)=0$ yields the condition $p_{22}^-=p_{22}^+$.
Then, the $3$rd generalized Lyapunov constant is given by
\begin{equation*}
\begin{aligned}
V_3(\varepsilon)=&-\frac{\pi}{8}
\big[3(a_3^++a_3^-)+(3p_{23}^++3p_{23}^--2p_{22}^+q_{22}^+
-2p_{22}^+q_{22}^-)\varepsilon^2 \big]\, \varepsilon^3.
\end{aligned}
\end{equation*}
Letting the $\varepsilon^{3}$-order and $\varepsilon^5$-order terms
in $V_3(\varepsilon)$ equal zero we obtain the conditions,
\begin{equation*}
\begin{aligned}
a_{3}^-=&-a_3^+,\\
p_{23}^-=&\ \frac{1}{3}\big(-3p_{23}^++2p_{22}^+q_{22}^+
+2p_{22}^+q_{22}^-\big).
\end{aligned}
\end{equation*}
Then, we have the 4th generalized Lyapunov constant, given by
\begin{equation*}
\begin{aligned}
V_4(\varepsilon)=& \frac{8}{45}\,
\big[ 3(b_3^++b_3^-)p_{22}^+ + 6a_3^+(q_{22}^++q_{22}^-) + M_{18}
\varepsilon^2\big]\, \varepsilon^5,
\end{aligned}
\end{equation*}where
\begin{equation*}
\begin{aligned}
M_{18}=6 p_{23}^+q_{22}^+-4p_{22}^+(q_{22}^+)^2+6p_{23}^+q_{22}^-
-4p_{22}^+q_{22}^+q_{22}^-+3p_{22}^+q_{23}^+-3p_{22}^+q_{23}^-.
\end{aligned}
\end{equation*}

Setting the $\varepsilon^5$-order term in
$V_4(\varepsilon)$ zero yields three necessary center conditions:
(i) $a_3^+=0$ and $b_3^-=b_3^+$, which leads to
conditions included in the condition VI;
(ii) $a_3^+=p_{22}^+=0$, giving the condition VII;
and (iii) $a_3^+\neq0$ and $q_{22}^-=\frac{1}{2a_3^+}(-b_3^+p_{22}^+
+ b_3^- p_{22}^+-2a_3^+q_{22}^+)$.

For the case (iii), by $M_{18}=0$ we have
\begin{equation*}
\begin{aligned}
q_{23}^+=\frac{1}{3a_3^+}(3b_3^+p_{23}^+-3b_3^-p_{23}^+
- 2 b_3^+p_{22}^+q_{22}^++ 2 b_3^- p_{22}^+q_{22}^++ 3 a_3^+q_{23}^-).
\end{aligned}
\end{equation*}
Then, we have the 5th generalized Lyapunov constant, given by
\begin{equation*}
\begin{aligned}
V_5(\varepsilon)=&-\frac{5\pi}{1152(a_3^+)^2}(b_3^+- b_3^-)\,
\big[54 (a_3^+)^3-a_3^+M_{19}\varepsilon^2-M_{20}\varepsilon^4\big]
\, \varepsilon^5,
\end{aligned}
\end{equation*}
where
\begin{equation*}
\begin{aligned}
M_{19}=&-7 b_3^+(p_{22}^+)^2 - b_3^- (p_{22}^+)^2
- 36 a_3^+p_{23}^++ 12 a_3^+p_{22}^+q_{22}^+,\\[0.5ex]
M_{20}=&-21 b_3^+(p_{22}^+)^2 p_{23}^++ 21 b_3^- (p_{22}^+)^2 p_{23}^+
- 54 a_3^+(p_{23}^+)^2 + 14 b_3^+ (p_{22}^+)^3 q_{22}^+\\
&-14 b_3^- (p_{22}^+)^3q_{22}^++36a_3^+p_{22}^+p_{23}^+q_{22}^+
-24a_3^+(p_{22}^+)^2 q_{23}^-.
\end{aligned}
\end{equation*}
Since $a_3^+\neq0$, we have $b_3^+=b_3^-$ from $V_5(\varepsilon)=0$.
Combining the conditions from the subcases (i) and (iii),
we have the condition VI.

\subsection{Case 4: $a_2^+a_2^-\neq0$}

Setting each $\varepsilon^k$-order term in
$V_2(\varepsilon)$ zero yields the necessary center
conditions $a_2^-=a_2^+$ and $p_{22}^-=p_{22}^+$.
Then, the $3$rd generalized Lyapunov constant is given by
\begin{equation*}
\begin{aligned}
V_3(\varepsilon)=&- \frac{\pi}{8}
\big[ 3a_3^+ + 3a_3^- - 2 a_2^+q_{22}^+ - 2 a_2^+q_{22}^-
+(3 p_{23}^++3 p_{23}^- - 2 p_{22}^+q_{22}^+-2 p_{22}^+q_{22}^-)\,
\varepsilon^2 \big]\, \varepsilon^3.
\end{aligned}
\end{equation*}
Letting the $\varepsilon^{3}$-order and $\varepsilon^5$-order terms
in $V_3(\varepsilon)$ equal zero we obtain the conditions,
\begin{equation*}
\begin{aligned}
a_{3}^-=&\ \frac{1}{3}(-3 a_3^++ 2 a_2^+q_{22}^+ + 2 a_2^+ q_{22}^-),\\
p_{23}^-=&\ \frac{1}{3}\big(3 p_{23}^+ + 2 p_{22}^+ q_{22}^+
+ 2 p_{22}^+ q_{22}^-\big).
\end{aligned}
\end{equation*}
Then, 4th generalized Lyapunov constant becomes
\begin{equation*}
\begin{aligned}
V_4(\varepsilon)=&\ \frac{1}{90}\,
\big[48a_2^+(b_3^++b_3^-)+16M_{21}\varepsilon^2
+225 a_2^+p_{23}^+\pi\varepsilon^3+16M_{22}\varepsilon^4
+225p_{22}^+p_{23}^+ \pi\varepsilon^5 \big]\, \varepsilon^3,
\end{aligned}
\end{equation*}where
\begin{equation*}
\begin{aligned}
M_{21}=&\ 3 b_3^+ p_{22}^+ - 3 b_3^- p_{22}^+ + 6 a_3^+ q_{22}^+
- 4 a_2^+ (q_{22}^+)^2 + 6 a_3^+ q_{22}^-- 4 a_2^+ q_{22}^+ q_{22}^-
+ 3 a_2^+ q_{23}^+ - 3 a_2^+ q_{23}^-,\\[0.5ex]
M_{22}=&\ 6 p_{23}^+ q_{22}^+ - 4 p_{22}^+ (q_{22}^+)^2
- 6 p_{23}^+ q_{22}^-- 4 p_{22}^+ q_{22}^+ q_{22}^-+ 3 p_{22}^+ q_{23}^+
- 3 p_{22}^+ q_{23}^-.
\end{aligned}
\end{equation*}
Letting the $\varepsilon^k$-order terms ($k=3,5,6,7$)
in $V_4(\varepsilon)$ equal zero we obtain the following conditions:
\begin{equation*}
\begin{aligned}
b_3^-=b_3^+, \quad p_{23}^+=0, \quad
q_{23}^-=\frac{1}{3a_2^+}\big[6 a_3^+ q_{22}^+ - 4 a_2^+ (q_{22}^+)^2 + 6 a_3^+ q_{22}^-- 4 a_2^+ q_{22}^+ q_{22}^-+ 3 a_2^+ q_{23}^+\big].
\end{aligned}
\end{equation*}
Then, $M_{22}$ is reduced to
$$
M_{22}=-\frac{6}{a_2^+}a_3^+ p_{22}^+ (q_{22}^+ + q_{22}^-).
$$
If $q_{22}^-=-q_{22}^+$, we obtain $a_3^-=-a_3^+$,
yielding the condition VIII.
Otherwise, we obtain the 5th  generalized Lyapunov constant, given by
\begin{equation*}
\begin{aligned}
V_5(\varepsilon)=&-\frac{1}{1440a_2^+}(q_{22}^+ + q_{22}^-)\,
\big[300 (a_2^+)^2 b_3^+ \pi+25\pi M_{23}\varepsilon^2
-6144 a_2^+ a_3^+ p_{22}^+\varepsilon^3\\
&-50p_{22}^+\pi M_{24}\varepsilon^4-6144 a_3^+ (p_{22}^+)^2 \pi\varepsilon^5
\big]\, \varepsilon^5,
\end{aligned}
\end{equation*}
where
\begin{equation*}
\begin{aligned}
M_{23}=&\ 27 (a_3^+)^2+12a_2^+b_3^+p_{22}^+-15a_2^+a_3^+q_{22}^+
-2(a_2^+)^2(q_{22}^+)^2+3a_2^+a_3^+q_{22}^-\\
&-2(a_2^+)^2q_{22}^+q_{22}^-+12(a_2^+)^2q_{23}^+,\\[0.5ex]
M_{24}=&\ 9 a_3^+q_{22}^++a_2^+(q_{22}^+)^2-12a_3^+q_{22}^-
+a_2^+q_{22}^+q_{22}^--6 a_2^+q_{23}^+.
\end{aligned}
\end{equation*}
With $b_3^+> \frac{1}{3}(a_2^+)^2$, we have $V_5(\varepsilon)\neq0$
when $a_3^+ p_{22}^+ =0$ and $q_{22}^+ + q_{22}^-\neq0$.

To this end, we finish the proof for the necessity of the conditions
in Theorem \ref{Th1}. Next, we prove the sufficiency of these eight
nilpotent center conditions.

If the condition I in Theorem \ref{Th1} holds, \eqref{Eqn1-5} is reduced to
{\begin{equation}\label{Eqn3-6}
(\dot{x},\, \dot{y})=\left\{\begin{aligned}
&\Big(y, \ -(b_2^+x^2+b_3^+x^3)\Big), \ \ &\text{if} \ \ x\geq0,\\[0.0ex]
&\Big(y, \ -b_3^-x^3\Big), \ \ &\text{if} \ \ x<0.
\end{aligned}\right.
\end{equation}
The two systems in \eqref{Eqn3-6}
are Hamiltonian systems, having respectively the Hamiltonian functions,
\begin{equation*}
\begin{aligned}
I^+(x,y)=\frac{1}{2}y^2+\frac{1}{3}b_2^+x^3+\frac{1}{4}b_3^+x^4
\quad \textrm{and} \quad
I^-(x,y)=\frac{1}{2}y^2+\frac{1}{4}b_3^-x^4,
\end{aligned}
\end{equation*}
which implies that the condition $I^+(0,y)\equiv I^-(0,y)$
in Proposition \ref{p11} is satisfied. So the origin
of system \eqref{Eqn3-6} is a nilpotent center.

Similar to the proof for the condition I, we can derive
the Hamiltonian functions (which are polynomials similar to the
above $I^{\pm}(x,y)$) for system \eqref{Eqn1-5} with
the condition II, or V or VII to
show that the origin of \eqref{Eqn1-5} is a nilpotent center
for these three cases.

If the condition III in Theorem \ref{Th1} holds,
\eqref{Eqn1-5} becomes
\begin{equation}\label{Eqn3-7}
(\dot{x},\, \dot{y})=\left\{\begin{aligned}
&\big(y-a_3^+x^3, \ -b_2^+x^2\big), \ \ &\text{if} \ \ x\geq0,\\[0.0ex]
&\big(y-a_3^-x^3, \ -b_2^-x^2\big), \ \ &\text{if} \ \ x<0,
\end{aligned}\right.
\end{equation}
with $-\,\frac{a_3^-}{a_3^+}=-\,\frac{b_2^-}{b_2^+} \equiv
h >0,\, \ne 1$. Then, system \eqref{Eqn3-7} can be rewritten as
\begin{equation}\label{Eqn3-7a}
(\dot{x},\, \dot{y})=\left\{\begin{aligned}
&\big(y-a_3^+x^3, \ -b_2^+x^2\big), \ \ &\text{if} \ \ x\geq0,\\[0.0ex]
&\big(y+ h\, a_3^+x^3, \ h\, b_2^+x^2\big), \ \ &\text{if} \ \ x<0.
\end{aligned}\right.
\end{equation}
Then, we apply the integrating factors:
$$
\mu^+ = \dfrac{9 (a_3^+)^2}{3 (a_3^+)^2 x^3-3 a_3^+ y-b_2^+}
\quad \textrm{and} \quad
\mu^- = \dfrac{9 (a_3^+)^2}{-3 (a_3^+)^2 h x^3-3 a_3^+ y-b_2^+}
$$
respectively to the first and second systems of \eqref{Eqn3-7a} to obtain
the first integrals:
$$
\left\{
\begin{array}{ll}
I^+(x,y) = -\,3 a_3^+ y + b_2^+ \ln|\! -\! 3 (a_3^+)^2 x^3+3 a_3^+ y+b_2^+|, &
\textrm{for} \ \ x \ge 0, \\[1.0ex]
I^-(x,y) = -\,3 a_3^+ y + b_2^+ \ln|3 h (a_3^+)^2 x^3+3 a_3^+ y+b_2^+|, &
\textrm{for} \ \ x < 0,
\end{array}
\right.
$$
which clearly shows that $ I^+(0,y) \equiv I^-(0,y)$, implying that
the origin of the system \eqref{Eqn3-7a} and so the system \eqref{Eqn3-7}
is a center.
%Note that the integrating factors $\mu^+$ and $\mu^-$
%keep the sign of the time $t$ unchanged.

If the condition IV in Theorem \ref{Th1}
holds, \eqref{Eqn1-5} can be rewritten as
\begin{equation}\label{Eqn3-8}
(\dot{x},\, \dot{y})=\left\{\begin{aligned}
&\Big(y-(a_2^+x^2+a_3^+x^3), \ -(b_2^+x^2+b_3^+x^3)\Big), \ \ &\text{if} \ \ x\geq0,\\[0.0ex]
&\Big(y-(a_2^+x^2-a_3^+x^3), \ -(-b_2^+x^2+b_3^+x^3)\Big), \ \ &\text{if} \ \ x<0,
\end{aligned}\right.
\end{equation}
which is symmetric with respect to the $y$-axis, and
by Proposition \ref{p22} we know that the origin
of \eqref{Eqn3-8} is a nilpotent center.
Similarly, since \eqref{Eqn1-5} is symmetric with respect to the $y$-axis
when the condition VI or VIII holds, the origin
of system \eqref{Eqn1-5} is a nilpotent center.

In the following, we present examples for the above systems
\eqref{Eqn3-6}, \eqref{Eqn3-7} and \eqref{Eqn3-8}
to show global phase portraits with a nilpotent center at the origin.

\begin{example}
The global phase portraits of the switching Li\'enard systems \eqref{Eqn3-6}, \eqref{Eqn3-7}
and \eqref{Eqn3-8} with three sets of parameter values
show a nilpotent center at the origin, as illustrated
in the Poincar\'e disc of Figure \ref{Fig2}.
\end{example}

\begin{figure}[!t]
\vspace{0.00in}
\centering
\hspace*{-0.4in}
\subfigure[]{
\label{Fig2a}
\includegraphics[width=0.28\textwidth,height=0.20\textheight]{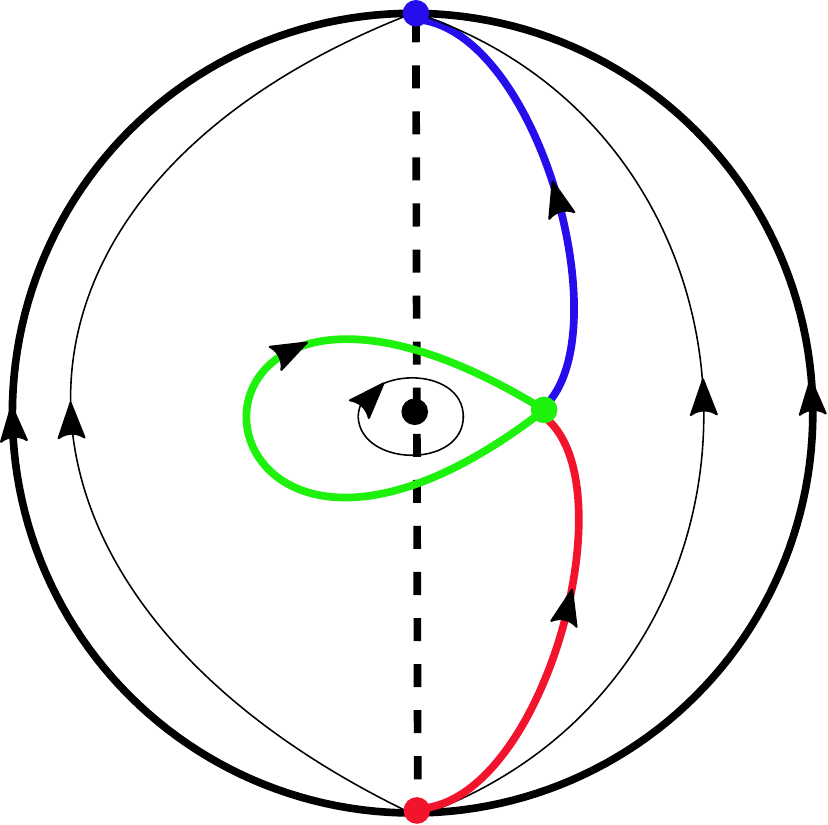}}
\hspace*{0.15in}
\subfigure[]{
\label{Fig2b}
\includegraphics[width=0.28\textwidth,height=0.20\textheight]{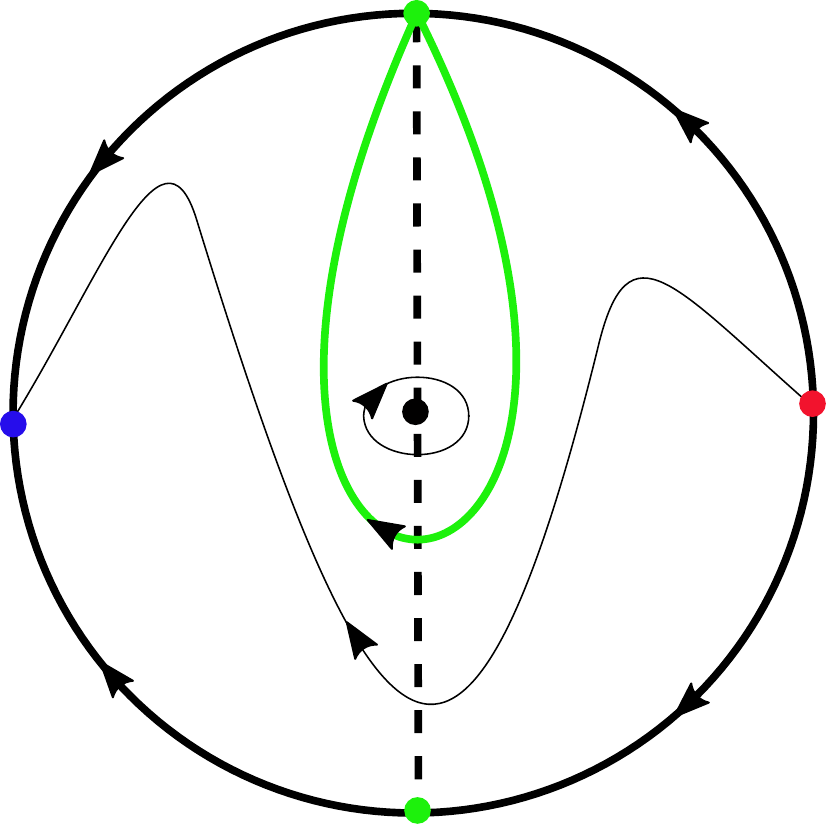}}
\hspace*{0.15in}
\subfigure[]{
\label{Fig2bc}
\includegraphics[width=0.28\textwidth,height=0.20\textheight]{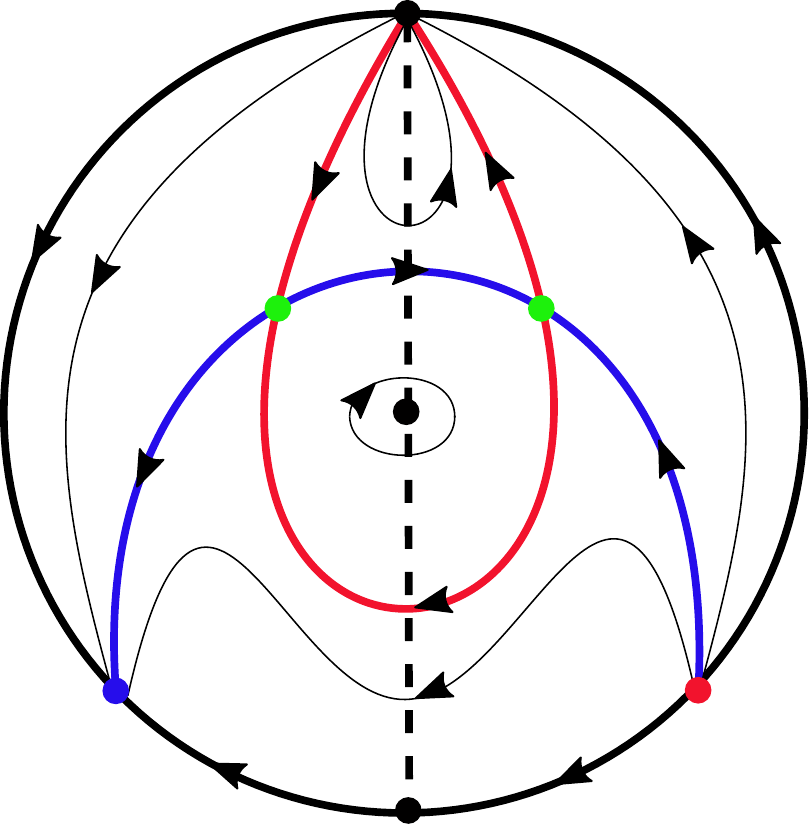}}
\hspace*{-0.3in}
\caption{Phase portraits for (a) system \eqref{Eqn3-6}
with $b_2^+=1,\,b_3^+=-2,\,b_3^-=1$,
(b) system \eqref{Eqn3-7} with $a_3^+=1,\, b_2^+=2,\, a_3^-=-2,\, b_2^-=-4$,
and (c) system \eqref{Eqn3-8} with $a_2^+=1,\, a_3^+=2,\, b_2^+=1,\,b_3^+=-2$.}
\label{Fig2}
\end{figure}

\section{The proof of Theorem \ref{Th2}}

We recall the notion of Poincar\'e compactification of switching differential
systems. This compactification identifies $\mathbb{R}^2$ with the interior
of the closed unit Poincar\'e disc $\mathbb{D}^2$ centered at the origin
of coordinates, and extends the switching differential systems to the
boundary of $\mathbb{D}^2$, which is called the circle of the infinity.  The switching manifold is indicated by
the dash line $``- - -"$ in $\mathbb{D}^2$. The singular points on the boundary of Poincar\'e disc are called {\it infinite
singular points}. More details about the Poincar\'e compactification for
switching polynomial systems can be found in \cite{Chen2022}.

With the results in \cite{LlibreV2022}, we have the
proposition which characterizes the global center of
switching polynomial systems.

\begin{proposition}\label{p33} {\rm [Propositions 5, \cite{LlibreV2022}]}
Consider the switching polynomial system \eqref{Eqn2-1}
which has a unique finite center at the singular point.
Assume that the infinity of \eqref{Eqn2-1} is not filled of singular points.
Then, the center is global if and only if all the infinite singular points
in $\mathbb{D}^2$, if they exist, are such that the local phase
portrait of each infinite singular point is consisting of two hyperbolic sectors, where the
two separatrices are on the circle of the infinity.
\end{proposition}

Suppose the origin of the switching Li\'enard system \eqref{Eqn1-5} is
the unique monodromic finite singular point.
Then, the remaining finite singular points, given by
\begin{equation*}
\begin{aligned}
e^+=(x^+,y^+)=\bigg(\!-\frac{b_2^+}{b_3^+},-\frac{(b_2^+)^2(a_3^+ b_2^+
- a_2^+ b_3^+)}{(b_3^+)^3}\bigg)
\end{aligned}
\end{equation*}and
\begin{equation*}
\begin{aligned}
e^-=(x^-,y^-)=\bigg(\! -\frac{b_2^-}{b_3^-},-\frac{(b_2^-)^2(a_3^- b_2^-
- a_2^- b_3^-)}{(b_3^-)^3}\bigg).
\end{aligned}
\end{equation*}
in the first and the second systems of \eqref{Eqn1-5} must be virtual,
i.e., $x^+<0$ and $x^->0$.
Hence, we have the necessity conditions:
\begin{equation}\label{Eqn4-1}
\begin{aligned}
G_1^+=&\ \{b_2^+=0,\ b_3^+\neq0\},\\[0.5ex]
G_2^+=&\ \{b_2^+>0,\ b_3^+\geq0\}
\end{aligned}
\end{equation}and
\begin{equation}\label{Eqn4-2}
\begin{aligned}
G_1^-=&\ \{b_2^-=0,\ b_3^-\neq0\},\\[0.5ex]
G_2^-=&\ \{b_2^-<0,\  b_3^-\geq0\}.
\end{aligned}
\end{equation}

Next, we consider the infinite singular point of \eqref{Eqn1-5}.
For studying the dynamical behaviours in the circle of the infinity,
we use the local charts: $U_{k}=\{(z_1,z_2)\in \mathbb{D}^2:z_k>0\}$
and $W_{k}=\{(z_1,z_2)\in \mathbb{D}^2:z_k<0\}$, $k=1,2$,
with the corresponding diffeomorphisms,
\begin{equation}
\phi_k:~U_k\rightarrow\mathbb{R}^2,\qquad \psi_k:~W_k
\rightarrow\mathbb{R}^2,
\end{equation}
defined by $\phi_k(z_1,z_2)=\big(\frac{z_2}{z_1},\frac{1}{z_1} \big)
=(u,w)$ and $\psi_k(z_1,z_2)=\big(\frac{z_1}{z_2},\frac{1}{z_2}\big)
=(u,w)$. Here, the coordinates $(u,w)$ play different roles in
the distinct local charts. Thus, the corresponding
vector fields of \eqref{Eqn1-5} in the local charts $U_1$ and
$W_1$ are given by
\begin{equation}\label{Eqn4-3}
\left(\begin{array}{cc}
\dot{u}\\
\dot{w}
\end{array}\right)
=\left(\begin{array}{c} -b_3^++a_3^+u-b_2^+w+a_2^+uw-u^2w^2 \\[0.5ex]
w(a_3^++a_2^+w-uw^2)
\end{array}\right)
\end{equation}
and
\begin{equation}\label{Eqn4-4}
\left(\begin{array}{cc}
\dot{u}\\
\dot{w}
\end{array}\right)
=\left(\begin{array}{c} -b_3^-+a_3^-u-b_2^-w+a_2^-uw-u^2w^2 \\[0.5ex]
w(a_3^-+a_2^-w-uw^2)
\end{array}\right),
\end{equation}
respectively. We derive the possible infinite singular points
$\big(\frac{b_3^{\pm}}{a_3^{\pm}},0 \big)$ or
$\big(\frac{b_2^{\pm}}{a_2^{\pm}},0 \big)$
which are nodes.
It follows from Proposition \ref{p33} that the origin of \eqref{Eqn1-5} cannot be a global
center in these cases. Hence, the systems \eqref{Eqn4-3} and \eqref{Eqn4-4}
cannot have these infinite singular points when the origin of \eqref{Eqn1-5} is a global center.
Therefore, the conditions \eqref{Eqn4-1} and \eqref{Eqn4-2} become
\begin{equation}\label{Eqn4-6}
\begin{aligned}
G_1^+=&\ \{a_2^+=a_3^+=b_3^+=0,\,\,b_2^+>0\},\\[0.5ex]
G_2^+=&\ \{a_3^+=0,\,\,b_2^+>0,\,\,b_3^+>0\},\\[0.5ex]
G_3^+=&\ \{a_3^+=0,\,\,b_2^+=0,\,\,b_3^+\neq0\}
\end{aligned}
\end{equation}and
\begin{equation}\label{Eqn4-7}
\begin{aligned}
G_1^-=&\ \{a_2^-=a_3^-=b_3^-=0,\,\,b_2^-<0\},\\[0.5ex]
G_2^-=&\ \{a_3^-=0,\,\,b_2^-<0,\,\,b_3^->0\},\\[0.5ex]
G_3^-=&\ \{a_3^-=0,\,\,b_2^-=0,\,\,b_3^-\neq0\}.
\end{aligned}
\end{equation}
Then we consider the following cases.

(i) Firstly, we consider the case $G_1^+\cap G_1^-$,
for which system \eqref{Eqn1-5} becomes
a quadratic switching Li\'enard
system. Since the expression of the quadratic switching
system \eqref{Eqn1-5} in $W_2$ is
the one in $U_2$ multiplied by $-1$,
the flows in the local phase portrait around the
origin of $W_2$ have the opposite direction compared
to the ones in $U_2$. Thus, we only need to study the
dynamics around the origin in $U_2$.
Then, in the local chart $U_2$ we have the following system,
\begin{equation}\label{Eqn4-8}
\left(\begin{array}{cc}
\dot{u}\\
\dot{w}
\end{array}\right)
=\left\{\begin{aligned}&\left(
\begin{array}{c} w+ b_2^+ u^3\\[0.5ex]
b_2^+ u^2 w
\end{array}\right), &\text{if} \ \ u\geq0,\\[1.0ex]
&\left(\begin{array}{c}
w+ b_2^- u^3\\[0.5ex]
b_2^- u^2 w
\end{array}\right), &\text{if} \ \ u<0,
\end{aligned}
\right.
\end{equation}
where both the two origins of the first and the second system
of \eqref{Eqn4-8} are nilpotent singular points.
Note that if the origin of $U_2$ is an infinite singular point,
it must come from the combination of the two origins of the first
and the second systems of \eqref{Eqn4-8}.
Since $b_2^+>0$ and $b_2^-<0$, it is easy to verify that
the origin of the first system of \eqref{Eqn4-8} is an unstable node
while the origin of the second one of \eqref{Eqn4-8} is a stable node.
Thus, it follows from Proposition \ref{p33} that
the origin of \eqref{Eqn1-5} cannot be a global center.

The above analysis shows that the infinite singular point
of the quadratic switching Li\'enard system \eqref{Eqn1-5} in $U_2$
does not consist of two hyperbolic sectors.
This implies that three is the lowest degree of the first and
the second systems of the switching Li\'enard system \eqref{Eqn1-5}
to have a global nilpotent center.
Thus, we only need to consider the case $G_i^+\cap G_j^-$, $i,j=2,3$.

(ii) For the case $G_i^+\cap G_j^-$, $i,j=2,3$,
we obtain the following system from \eqref{Eqn1-5},
\begin{equation}\label{Eqn4-9}
\left(\begin{array}{cc}\dot{u}\\
\dot{w}
\end{array}\right)
=\left\{\begin{aligned}&\left(
\begin{array}{c} b_3^+u^4-a_2^+u^2w+b_2^+u^3w+w^2 \\[0.5ex]
u^2w(b_3^+u+b_2^+w)
\end{array}\right), &\text{if} \ \ u\geq0,\\[1.0ex]
&\left(\begin{array}{c}
b_3^-u^4-a_2^-u^2w+b_2^-u^3w+w^2 \\[0.5ex]
u^2w(b_3^-u+b_2^- w)
\end{array}\right), &\text{if} \ \ u<0,
\end{aligned}
\right.
\end{equation}
in $U_2$ and $W_2$. Since system \eqref{Eqn1-5} is a cubic switching system, the local qualitative property of the
origin of $W_2$ has the same sense with respect
to the one in $U_2$. Obviously, the origins of the first and
second system in \eqref{Eqn4-9} in $U_2$ are singular points
which are linearly zero.
In order to understand the local qualitative properties around
these two singular points, we apply a direct blow-up
$(u,v)\rightarrow(u,w)$ with $v= \frac{w}{u}$ to eliminate
the common factor $u$, yielding
\begin{equation}\label{Eqn4-10}
\left(\begin{array}{cc}\dot{u}\\
\dot{v}
\end{array}\right)
=\left\{\begin{aligned}&\left(\begin{array}{c}
u(b_3^+ u^2 - a_2^+  u v + b_2^+  u^2 v + v^2) \\[0.5ex]
(a_2^+  u - v) v^2
\end{array}\right), &\text{if} \ \ u\geq0,\\[1.0ex]
&\left(\begin{array}{c}
u(b_3^- u^2 - a_2^- u v + b_2^- u^2 v + v^2) \\[0.5ex]
(a_2^- u - v) v^2
\end{array}\right), &\text{if} \ \ u<0,
\end{aligned}
\right. \vspace{0.10in}
\end{equation}
where the two origins are singular points with linearly zero.
We do a further blow-up $(u,v)\rightarrow(u,V)$ with $V= \frac{v}{u}$,
to obtain the following system,
\begin{equation}\label{Eqn4-11}
\left(\begin{array}{cc}\dot{u}\\
\dot{V}
\end{array}\right)
=\left\{\begin{aligned}&\left(\begin{array}{c}
u(b_3^+ - a_2^+ V + b_2^+ u V + V^2) \\[0.5ex]
-V(b_3^+ - 2 a_2^+ V + b_2^+ u V + 2 V^2)
\end{array}\right), &\text{if} \ \ u\geq0,\\[1.0ex]
&\left(\begin{array}{c}
u(b_3^- - a_2^-  V + b_2^-  u V + V^2) \\[0.5ex]
-V(b_3^-  - 2 a_2^-  V + b_2^-  u V + 2 V^2)
\end{array}\right), &\text{if} \ \ u<0,
\end{aligned}
\right.
\end{equation}
where the common factor $u^2$ has been eliminated.

When $(a_2^+)^2-2b_3^+<0$, for $u=0$ the first system of \eqref{Eqn4-11}
has a saddle at $E_1=(0,0)$. Going back
through the change of variables we get that the local phase portrait
of the right origin in $U_2$  has two hyperbolic sectors, see Figure \ref{Fig3}.
When $(a_2^+)^2=2b_3^+\neq0$, for $u=0$ the first system of \eqref{Eqn4-11}
has two singular points at $E_1=(0,0)$ (a saddle)
and at $E_2=(\frac{a_2^+}{2},0)$ (a saddle-node). Again, going back
through the change of variables we find that the local phase portrait of
the right origin in $U_2$ has two hyperbolic and two parabolic
sectors, see Figure \ref{Fig4}.
When $(a_2^+)^2-2b_3^+>0$, for $u=0$ the first system of \eqref{Eqn4-11}
has three singular points at $E_1=(0,0)$ (a saddle),
$E_{2,3}=\big(\frac{1}{2} (a_2^+\pm ((a_2^+)^2-2b_3^+)^{1/2},0 \big)$
(a saddle and a node). Going back
through the change of variables we obtain that the local qualitative property
of the right origin in $U_2$ has two hyperbolic and
two parabolic sectors, see Figure \ref{Fig5}.

\begin{figure*}[!h]
\begin{center}
\begin{overpic}[width=0.80\textwidth]{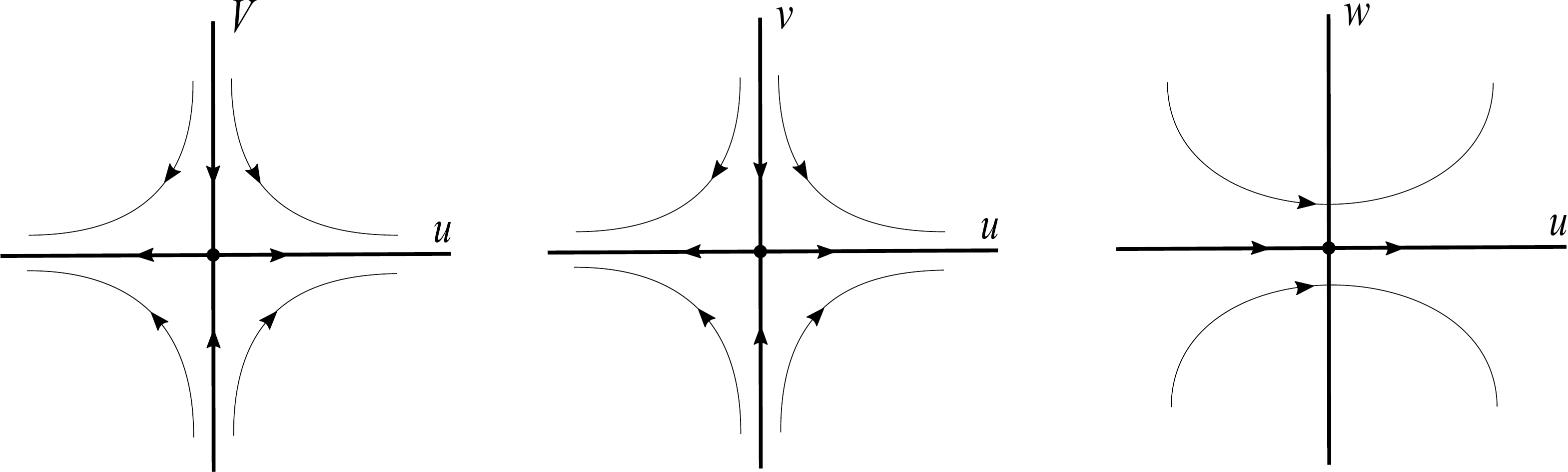}
\put(11,-5){$(a)$}
\put(47,-5){$(b)$}
\put(83,-5){$(c)$}
\end{overpic}
\vspace{3mm}
\end{center}
\caption{Local phase portraits around the origin of the local
chart $U_2$ with $(a_2^+)^2-2b_3^+<0$ for
(a) the first system of  \eqref{Eqn4-11}, (b) the first system of  \eqref{Eqn4-10}, and
(c) the first system of  \eqref{Eqn4-9}.}
\label{Fig3}
\end{figure*}

\begin{figure*}[!h]
\begin{center}
\begin{overpic}[width=0.80\textwidth]{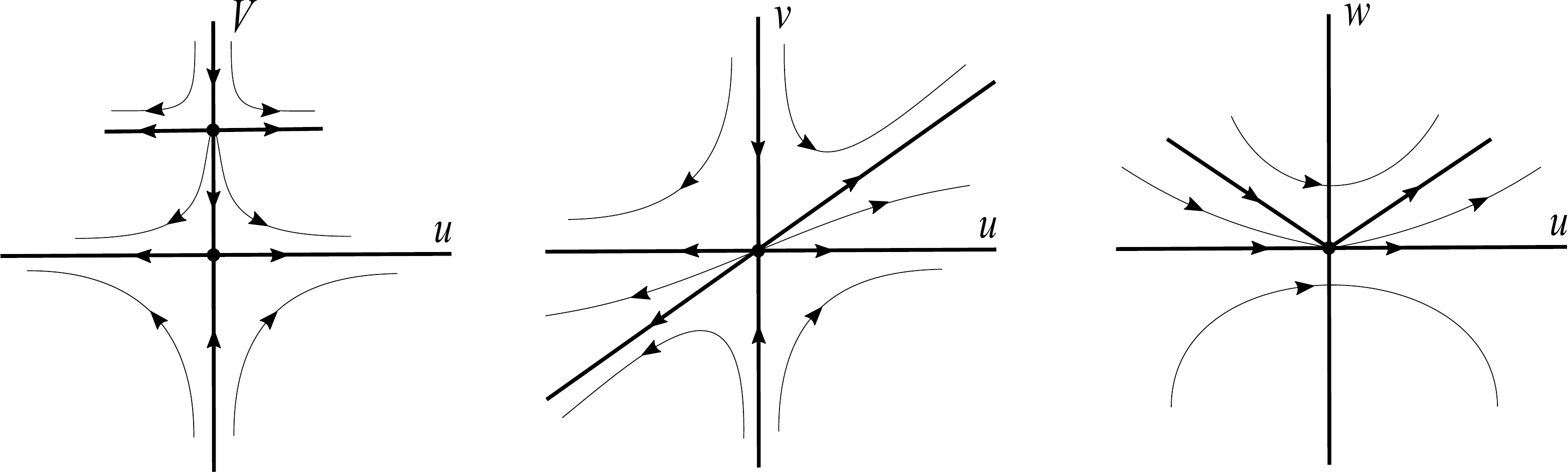}
\put(11,-5){$(a)$}
\put(47,-5){$(b)$}
\put(83,-5){$(c)$}
\end{overpic}
\vspace{3mm}
\end{center}
\caption{Local phase portraits around the origin of the local chart
$U_2$ with $(a_2^+)^2-2b_3^+=0$ for
(a) the first system of  \eqref{Eqn4-11}, (b) the first system of  \eqref{Eqn4-10}, and
(c) the first system of  \eqref{Eqn4-9}.}
\label{Fig4}
% Give a unique label
\end{figure*}

\begin{figure*}[!h]
\begin{center}
\begin{overpic}[width=0.80\textwidth]{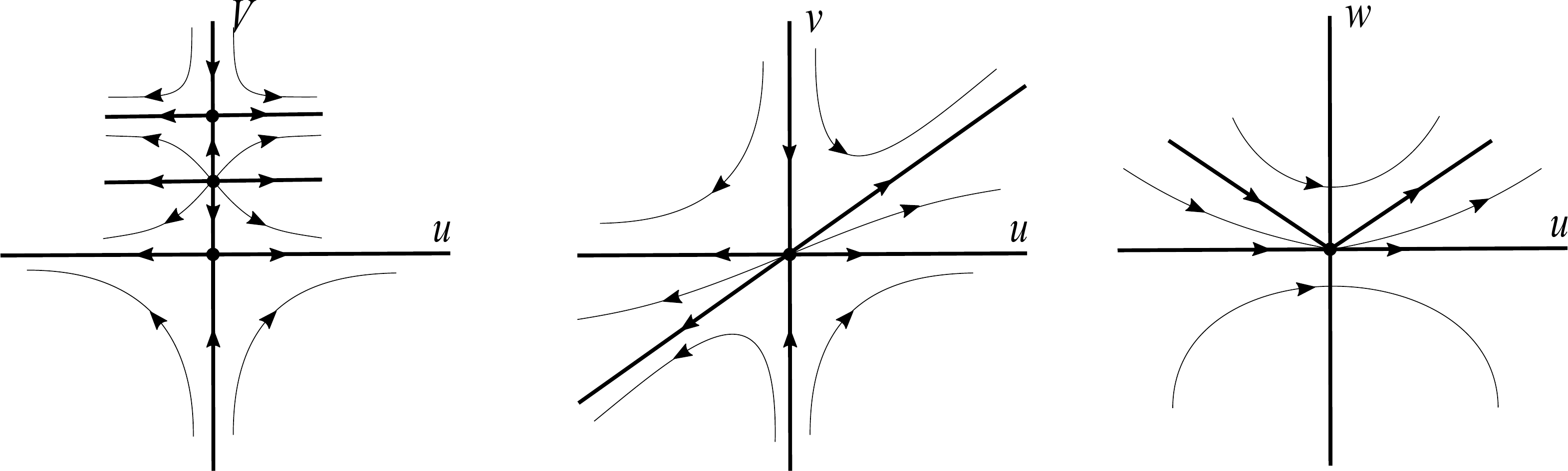}
\put(11,-5){$(a)$}
\put(47,-5){$(b)$}
\put(83,-5){$(c)$}
\end{overpic}
\vspace{3mm}
\end{center}
\caption{Local phase portraits at the origin of the local chart $U_2$
 with $(a_2^+)^2-2b_3^+>0$ for
(a) the first system of \eqref{Eqn4-11}, (b) the first system \eqref{Eqn4-10},
and (c) the first system of \eqref{Eqn4-9}.}
\label{Fig5}
% Give a unique label
\end{figure*}

Similarly, we derive that the local qualitative property around the left
origin in $W_2$ has two hyperbolic sectors if and only if
$(a_2^-)^2-2b_3^-<0$. Thus, the switching Li\'enard system \eqref{Eqn1-5}
can only have a global center at the origin if condition
\begin{equation}\label{Eqn4-12}
\begin{aligned}
G=&\{a_3^{\pm}=0,\ b_2^+\geq0,\ b_2^-\leq0,\ (a_2^+)^2-2b_3^+<0,\
(a_2^-)^2-2b_3^-<0\}
\end{aligned}
\end{equation}
holds. Thus, by Theorem \ref{Th1} we have proved Theorem \ref{Th2}.

In fact, combining the conditions I, II, V and VII of Theorem \ref{Th1}
and the above condition $G$, we obtain the condition $G_1$ of
Theorem \ref{Th2}. From the condition III of Theorem \ref{Th1},
we have $a_3^+a_3^-\neq0$ which contradicts the condition $G$, and thus
the origin of system \eqref{Eqn1-5} cannot be a global center.
From the conditions IV, VI and VIII of Theorem \ref{Th1} together with
the condition $G$, we obtain the condition $G_2$ of Theorem \ref{Th2}.

\section{The proof of Theorem \ref{Th3}}

In this section, we prove the existence of $5$ small-amplitude limit
cycles by perturbing the nilpotent origin of the cubic switching Li\'enard system \eqref{Eqn1-5}
under the center condition IV.
More precisely, with the center condition IV, we perturb system
\eqref{Eqn1-5} to obtain
\begin{equation}\label{Eqn5-1}
\left(\begin{array}{cc}\dot{x}\\
\dot{y}
\end{array}\right)
=\left\{\begin{aligned}&\left(\begin{array}{c}
y -(a_2^+x^2+a_3^+x^3)+d\varepsilon^{9}
+\displaystyle\sum_{k=1}^{6}\varepsilon^{k+1}
\delta_{k}x+\sum_{k=1}^{5}\varepsilon^k (p_{k2}^+x^2+p_{k3}^+x^3) \\[1.0ex]
-\varepsilon^2x-(b_2^+x^2+b_3^+x^3)
+\displaystyle\sum_{k=1}^{5}\varepsilon^k
(q_{k2}^+x^2+q_{k3}^+x^3)
\end{array}\right), &\text{if} \ \ x\geq0,\\[1.0ex]
&\left(\begin{array}{c}
y -(-a_2^+x^2+a_3^+x^3)
+\displaystyle\sum_{k=1}^{6}\varepsilon^{k+1} \delta_{k}x
+\sum_{k=1}^{5}\varepsilon^k (p_{k2}^-x^2+p_{k3}^-x^3) \\[1.0ex]
-\varepsilon^2x-(b_2^+x^2-b_3^+x^3)+\displaystyle\sum_{k=1}^{5}\varepsilon^k
(q_{k2}^-x^2+q_{k3}^-x^3)
\end{array}\right), &\text{if} \ \ x<0.
\end{aligned}
\right.
\end{equation}
Denote by one parameter vector $$\lambda^*=(a_2^+,a_3^+,b_2^+,b_3^+,d,\delta_i,p_{j2}^{\pm},p_{j3}^{\pm},q_{j2}^{\pm},q_{j3}^{\pm}),$$
$i=1,2,...,6$, $j=1,2,...,5$. Further, introducing the scaling
$(x,y,t)\rightarrow \big(\varepsilon^2 y, \varepsilon^3 x,
\frac{t}{\varepsilon} \big)$ into \eqref{Eqn5-1}, we obtain
a system up to $\varepsilon^{6}$-order terms,
\begin{equation}\label{Eqn5-2}
\!\!\! \left(\! \begin{array}{cc}
\dot{x}\\
 \dot{y}
\end{array} \! \right)
\!=\! \left\{
\begin{aligned}
&\left(\begin{array}{c}
-y-b_2^+y^2-b_3^+\varepsilon^2y^3-\displaystyle\sum_{k=1}^{5}
\varepsilon^k q_{k2}^+ y^2-\displaystyle\sum_{k=1}^{5} \varepsilon^{k+2}
q_{k3}^+ y^3 \\[1.0ex]
x+\displaystyle\sum_{k=1}^{6}\varepsilon^k \delta_{k}y-(a_2^+\varepsilon y^2
+a_3^+\varepsilon^3 y^3)-\displaystyle\sum_{k=1}^{5} \varepsilon^{k+1}
p_{k2}^+ y^2 -\displaystyle\sum_{k=1}^{5} \varepsilon^{k+3} p_{k3}^+ y^3
\end{array}\right)\!, \!\! & \text{if} \ y\geq0,\\
&\left(\!\!\!
\begin{array}{c}
-y+b_2^+y^2-b_3^+\varepsilon^2y^3-\displaystyle\sum_{k=1}^{5} \varepsilon^k
q_{k2}^- y^2-\displaystyle\sum_{k=1}^{5} \varepsilon^{k+2} q_{k3}^- y^3
\\[1.0ex]
x+ \! \displaystyle\sum_{k=1}^{6}\varepsilon^k \delta_{k}y
\!-\! (a_2^+\varepsilon y^2 \!-\! a_3^+\varepsilon^3 y^3)
\!+\! d \varepsilon^6 \!-\! \displaystyle\sum_{k=1}^{5}
\varepsilon^{k+1} p_{k2}^- y^2
\!-\! \displaystyle\sum_{k=1}^{5} \varepsilon^{k+3} p_{k3}^- y^3
\end{array} \!\!\! \right)\!, &\text{if} \ y<0.
\end{aligned}
\right.
\end{equation}
Note that we let the constant term be zero, i.e., $d=0$,
when we compute $\varepsilon^k$-order Lyapunov constants.
By using the generalized Lyapunov constants, we start from
the $\varepsilon$-order terms to determine the number of limit cycles.
Thus, we need to find $\varepsilon^k$-order center conditions such that all the
$\varepsilon^k$-order Lyapunov constants being zero when we want to
obtain more limit cycles by using $\varepsilon^{k+1}$-order terms.
That is to say that the origin of the system is a center up to
$\varepsilon^k$-order, see \cite{TY2018}.
Finally, we perturb the coefficient $d$ of the constant term
to obtain one more small-amplitude limit cycle.

More precisely, to prove the existence of small-amplitude limit cycles
around the origin, we compute the higher-order Lyapunov constants
$V_{jk}$ for $j=1,2,\ldots, 5$ and $ k=1,2,\ldots,6$.  We start from $k=1$.
For a fixed $k$, we choose appropriate parameter
values such that as many as possible higher-order Lyapunov constants vanish.

First, for $k=1$, we obtain that all $\varepsilon$-order Lyapunov
constants are zero when $\delta_1=0$, i.e.,
$V_{j1}=0$ for all $j$.

Next, for $k=2$, we have $V_{12}=\delta_2 \pi=0$ when $\delta_2=0$.
Then, we solve $V_{22}=\frac{4}{3} (p_{12}^- - p_{12}^+)=0$ to obtain
$$
p_{12}^-= p_{12}^+,
$$
under which the 3rd $\varepsilon^2$-order Lyapunov constant becomes
$$
V_{32}= \, \frac{\pi}{4}\,a_2^+ (q_{12}^- + q_{12}^+).
$$
Setting $q_{12}^- =- q_{12}^+$ yields $V_{j2}=0$ for all $j$.
So, for $k=2$, we can have $2$ limit cycles
by perturbing $\delta_2$ and $p_{12}^-$ with $V_{32} \ne 0$.

Now, consider $k=3$. We have $V_{13}=0$ with $\delta_3=0$ and
$V_{23}=\frac{4}{3} (p_{22}^- - p_{22}^+)$. Solving $V_{23}=0$ we have
$$
p_{22}^-= p_{22}^+.
\vspace{-0.10in}
$$
Further, the 3rd $\varepsilon^3$-order Lyapunov constant is obtained as
$$
V_{33}= \, \frac{\pi}{4}\,a_2^+ (q_{22}^- + q_{22}^+).
$$
Setting $q_{22}^- =- q_{22}^+$ yields $V_{j3}=0$ for all $j$,
which indicates that
for $k=3$
we can have $2$ limit cycles if choosing $V_{33} \ne 0$ and
perturbing $\delta_3$ and $p_{22}^-$.

For $k=4$, solving $V_{24}=V_{34}=0$ we have
$$
p_{32}^-= p_{32}^+ \quad \textrm{and} \quad
p_{13}^-=\frac{2}{3}a_2^+q_{32}^- +\frac{2}{3}a_2^+q_{32}^+- p_{13}^+,
$$
under which the 4th $\varepsilon^4$-order Lyapunov constant becomes
$$
V_{44}= -\, \frac{8}{45}\, a_2^+(4q_{32}^-b_2^+
+ 4 b_2^+q_{32}^+ + 3q_{13}^- - 3q_{13}^+).
$$
Setting $V_{44} = 0 $ we have
$$
q_{13}^- =-\frac{4}{3}b_2^+q_{32}^- -\frac{4}{3}b_2^+q_{32}^+- q_{13}^+,
$$
leading to $V_{j4}=0$ for all $j$. So, for $k=4$ we can have
$3$ limit cycles by choosing $V_{44} \ne 0 $ and
perturbing $\delta_4$, $p_{32}^-$ and $\ p_{13}^-$.

For $k=5$, we similarly obtain
$$
p_{42}^-= p_{42}^+
$$
by solving $V_{25}=0$, and
$$
p_{23}^-=\frac{2}{3}q_{32}^-p_{12}^+ + \frac{2}{3}q_{42}^-a_2^+
+ \frac{2}{3}a_2^+q_{42}^++\frac{2}{3}p_{12}^+q_{32}^+-p_{23}^+
$$
by solving $V_{35}=0$. Then, the 4th $\varepsilon^5$-order Lyapunov
constant is given by
$$
V_{45}=-\frac{8}{45}a_2^+(4q_{32}^-q_{12}^+ + 4q_{42}^-b_2^+ + 4b_2^+q_{42}^+ + 4q_{12}^+q_{32}^+ + 3Q_{23}^- - 3q_{23}^+).$$
Setting
$$
q_{23}^-=-\frac{4}{3}q_{12}^+q_{32}^- -\frac{4}{3}q_{42}^-b_2^+
-\frac{4}{3}b_2^+q_{42}^+ -\frac{4}{3}q_{12}^+q_{32}^+ + q_{23}^+,
$$
under which $V_{j5}=0$ for all $j$. This implies that
$3$ limit cycles can be obtained for $k=4$
if choosing $V_{45} \ne 0 $ and
perturbing $\delta_5$, $p_{42}^-$ and $p_{23}^-$.

Finally, for $k=6$, similarly by solving $V_{16}=V_{26}=V_{36}=0$ we obtain
$$
\begin{array}{ll}
\delta_6=0,\\[1.0ex]
p_{52}^-= p_{52}^+,\\[1.0ex]
p_{33}^-=\frac{2}{3}q_{32}^-p_{22}^+ +\frac{2}{3}q_{42}^-p_{12}^+
+ \frac{2}{3}q_{52}^-a_2^+ + \frac{2}{3}a_2^+q_{52}^+
+\frac{2}{3}p_{12}^+q_{42}^++\frac{2}{3}p_{22}^+q_{32}^+-p_{33}^+.
\end{array}
$$
Then, we solve $V_{46}=0$ for $q_{33}^-$ to obtain
\begin{equation*}
\begin{aligned}
q_{33}^-=&-\frac{4}{3a_2^+}(q_{32}^-a_2^+q_{22}^+ + 4q_{42}^-a_2^+q_{12}^+
+ 4q_{52}^-a_2^+b_2^+ + 4a_2^+b_2^+q_{52}^+ \\
&+ 4a_2^+q_{12}^+q_{42}^+ + 4a_2^+q_{22}^+q_{32}^+ - 6q_{32}^-a_3^+
- 3a_2^+q_{33}^+ - 6a_3^+q_{32}^+),
\end{aligned}
\end{equation*}
under which the 5th $\varepsilon^6$-order Lyapunov constant is reduced to
$$
V_{56} = -\dfrac{5 \pi}{48} \,
(2a_2^+b_3^+ - 3a_3^+b_2^+)(q_{32}^- + q_{32}^+).
$$
So letting
$$
(2a_2^+b_3^+ - 3a_3^+b_2^+)(q_{32}^- + q_{32}^+)=0
$$
yields $V_{j6}=0$ for all $j$.
Therefore, for $k=6$, we can obtain
$4$ limit cycles
for choosing $V_{56} \ne 0 $.

Summarizing the above results, $4$ small-amplitude
limit cycles can be obtained from the $\varepsilon^{6}$-order
Lyapunov constants. Assume that $p_{22}^+\neq0$,
a direct computation shows that
$$
\begin{array}{rl}
& \det \left[ \dfrac{\partial (V_{16},V_{26},V_{36},V_{46})}{\partial(\delta_6,p_{52}^-,q_{32}^+,q_{33}^+)}
\right]
= \dfrac{8\pi^2 }{45}a_2^+ p_{22}^+\ne 0,
\end{array}
$$
which implies that there exists $4$ small-amplitude limit cycles
by perturbing the parameters, in backward, $\delta_6$,
$p_{52}^-$, $ q_{32}^+$ and $q_{33}^+$.

One more small-amplitude limit cycle is obtained by perturbing $d$,
leading to a total $5$ small-amplitude limit cycles.
For $d$ and $\varepsilon$ sufficiently small, the switching
system \eqref{Eqn5-2} has a small sliding segment on the switching
manifold $y=0$ with the end points at
$(0,0)$ and $(x_\varepsilon,0)$,
where $x_\varepsilon=-d\varepsilon^6$ is the unique root of the
equation $\dot{y}=0$ in the second system of \eqref{Eqn5-2}.
Thus, the sliding segment shrinks to $(0,0)$ when $\varepsilon$ goes to zero.

For the point $(x_e,0)$ on the switching manifold $y=0$ with
$x_\varepsilon<x_e\ll1$, we define a bifurcation function,
\begin{equation}\label{Eqn5-3}
d(x_e,\varepsilon)=\Pi^+(x_e,\varepsilon)-\Pi_-^+(x_e,\varepsilon),
\end{equation}
for which has two small positive constants $\varepsilon_1$
and $\varepsilon_2$ such that $\Pi^+(x_e,\varepsilon):
(0,\varepsilon_1)\rightarrow(-\varepsilon_1,0)$ is the first half-bifurcation
function of \eqref{Eqn5-2},
and $\Pi_-^+(x_\varepsilon,\varepsilon):
(x_\varepsilon,x_\varepsilon+\varepsilon_2)
\rightarrow(x_\varepsilon-\varepsilon_2,x_\varepsilon)$ is the
second half-bifurcation function derived from the second system
of \eqref{Eqn5-2} by the transformation $(x,y,t)\rightarrow(x,-y,-t)$.
By the polar coordinate transformation $x=r\cos\theta$ and
$y=r\sin\theta$, we have
\begin{equation}\label{Eqn5-4}
\Pi^+(x_e,\varepsilon)=V_1^+(\varepsilon) x_e+O(x_e^2) \quad \text{and} \quad
\Pi_-^+(x_e,\varepsilon)=V_0^-(\varepsilon)+V_1^-(\varepsilon) x_e+O(x_e^2),
\end{equation}
where
\begin{equation}\label{Eqn5-5}
V_0^-(\varepsilon)=2d\varepsilon^6+O(\varepsilon^7) \quad \text{and} \quad
V_1^{\pm}(\varepsilon)=\pm\frac{1}{2}\delta_6\pi\varepsilon^6+O(\varepsilon^7),
\end{equation}
It follows from \eqref{Eqn5-3} and \eqref{Eqn5-4} that
$V_{06}=-2d$ and $V_{16}|_{d=0}=\delta_6\pi$. Hence, we obtain
the $\varepsilon^{6}$-order Lyapunov
constants $V_{06}$, $V_{16}$, $V_{26}$, $V_{36}$ and $V_{46}$, which
are independent of the parameter vector $\lambda^*$.
Therefore, we have proved the existence of $5$ small-amplitude
crossing limit cycles near the origin of cubic switching Li\'enard system \eqref{Eqn1-5}.

To end this section, we provide an example with exact critical parameter
vector $\lambda^*$ to demonstrate
the existence of $5$ small-amplitude limit cycles. To achieve this,
we need to find $5$ positive roots which are solved from the
displacement equation \eqref{Eqn2-9}:
\begin{equation}\label{Eqn6-1}
d(\xi)=V_0(\lambda^*)+V_1(\lambda^*)\xi + \cdots
+V_4(\lambda^*)\xi^4 + V_5(\lambda^*)\xi^5 +O(\xi^6)=0.
\end{equation}
As discussed above, we have
$V_j(\lambda^*)=V_{j6}\varepsilon^6+O(\varepsilon^6)$,
$j=0,1, \ldots,5$.
Hence, \eqref{Eqn6-1} becomes
\begin{equation}\label{Eqn6-2}
d(\xi)=(V_{06}+V_{16}\xi + \cdots +V_{46}\xi^4
+V_{56}\xi^5)\varepsilon^6+O(\varepsilon^7)=0.
\end{equation}
In general, it is a very challenging task to find the perturbed
parameter vector to give a numerical realization.
However, for our system the parameters in the displacement equation are
linearly decoupled. Thus, we can obtain these
values by perturbing exactly one parameter for one higher-order
Lyapunov constant
such that equation \eqref{Eqn6-1} has $5$ positive roots.

In the following, we construct a concrete example.
First, note that the critical parameter vector $\lambda^*$
given above satisfies that
$V_{j6}=0$, $j=0,1,\ldots,4$, but $V_{56}\neq0$, and $V_{56}$ does not contain any non-zero
lower-order $\varepsilon$ terms.
Thus, we can take perturbations in the backward order:
on $q_{32}^-$ for $V_{56}$, on $q_{33}^-$ for $V_{46}$,
on $p_{33}^-$ for $V_{36}$,
on $p_{52}^-$ for $V_{26}$, on $\delta_{6}$ for $V_{16}$,
and on $d$ for $V_{06}$.
In particular, with the help of Maple taking the accuracy up to
$60$ decimal points,
we set the free parameters to take the values:
\begin{equation}\label{Eqn6-3}
\begin{array}{llllllll}
&a_2^+=b_3^+=p_{12}^+=p_{22}^+=p_{33}^+=p_{52}^+=q_{12}^+=q_{22}^+
=q_{32}^+=q_{33}^+=q_{42}^{\pm}=q_{52}^{\pm}=1,\\[1.0ex]
&a_3^+-3=b_2^++2=q_{32}^-+2=0,
\end{array}
\end{equation}
and choose the following values for the perturbed parameters:
\begin{equation}\label{Eqn6-4}
\begin{array}{rlrl}
d^* \, =\!\!\!\!& d-5\times 10^{-31}, \, &\delta_6^* \ =\!\!\!\!&
\delta_6-3.2\times 10^{-21}, \\[0.5ex]
p_{52}^{-*} = \!\!\!\! & p_{52}^{+}+1.5\times 10^{-13},
\hspace*{0.20in} &p_{33}^{-*}
= \!\!\!\! & p_{33}^{-}-8.5\times 10^{-8},\\[0.5ex]
q_{33}^{-*} = \!\!\!\! & q_{33}^{-}+0.0188.
\end{array}
\end{equation}
With the above perturbed parameter values, we obtain the
perturbed $\varepsilon^{6}$-order Lyapunov constants as follows:
\begin{equation}\label{Eqn6-5}
\begin{array}{llllllll}
&V_{06}=-1.0 \times 10^{-30}, \quad &
V_{16}\approx 1.0053096491\times10^{-20}, \\[0.5ex]
&V_{26}\approx -2.0006218904\times10^{-13},\hspace*{0.20in}
&V_{36}\approx 1.0013826592\times10^{-7},\\[0.5ex]
&V_{46}\approx -0.0100266667, \quad &V_{56}\approx6.5449846949,
\end{array}
\end{equation}
for which the displacement equation \eqref{Eqn6-2} has $5$ positive roots:
\begin{equation}\label{Eqn6-6}
\!\!\! \begin{array}{llllllll}
\xi_1 \!\!\!& \approx 0.0015219219,
\quad &\xi_2 \!\!\!& \approx 7.3120695527\times 10^{-6},
\quad &\xi_3 \!\!\!& \approx 2.6762403416\times 10^{-6},
\\[1.0ex]
\xi_4 \!\!\!& \approx 5.1472016794\times 10^{-8},
\quad &\xi_5 \!\!\!& \approx 9.9669521943\times 10^{-11},
\end{array}
\end{equation}
which are approximations of the amplitudes for the bifurcating
$5$ small-amplitude limit cycles.

\section{Conclusion}

In this paper, we have studied the nilpotent center problem and
the limit cycle bifurcation problem for switching
nilpotent systems in $\mathbb{R}^2$.
We have developed a higher-order Poincar\'e-Lyapunov method to compute
the generalized Lyapunov constants for switching nilpotent systems.
By using this method, we derive the nilpotent center conditions for
the origin of the switching cubic polynomial Li\'enard systems.
Further, we characterize all the switching cubic polynomial
Li\'enard systems to have a global nilpotent center.
Finally, we construct a perturbed system with one of the center
conditions to show the existence of $5$ small-amplitude limit cycles
bifurcating from the nilpotent center, which is a new lower bound
of the maximal number of limit cycles in such switching cubic
Li\'enard systems with a nilpotent singular point.
The methodology developed in this paper
can be applied to investigate complex dynamics of other nonlinear
systems with nilpotent singular points.

\section*{Acknowledgment}

This work was supported by the National Natural Science
Foundation of China, No.~12001112 (T. Chen) and No.~12071198 (F. Li),
Guangdong Basic and Applied Basic Research Foundation,
No.~2022A1515011827 (T. Chen), and the Natural Sciences and Engineering
Research Council of Canada, No.~R2686A02 (P. Yu).


\begin{thebibliography}{99}

\bibitem{Adriana}
A. Buic$\rm{\breve{a}}$, J. Llibre, O. Makarenkov,
Asymptotic stability of periodic solutions for nonsmooth differential
equations with application to the nonsmooth van der Pol oscillator,
SIAM J. Math. Anal. 40 (2009) 2478--2495.

\bibitem{Caubergh}
M. Caubergh,
Hilbert's sixteenth problem for polynomial Li\'enard equations,
Qual. Theory Dyn. Syst. 11 (2012) 3--18.

\bibitem{CHB2015}
H. Chen, X. Chen, Dynamical analysis of a cubic Li\'enard system with global parameters,
Nonlinearity 28 (2015) 3535--3562.

\bibitem{CHB2016}
H. Chen, X. Chen, Dynamical analysis of a cubic Li\'enard system with global parameters (II),
Nonlinearity 29 (2016) 1798--1826.

\bibitem{CHB2020}
H. Chen, X. Chen, Dynamical analysis of a cubic Li\'enard system with global parameters III,
Nonlinearity 33 (2020) 1443--1465.

\bibitem{CDT2018}
H. Chen, S. Duan, Y. Tang, J. Xie,
Global dynamics of a mechanical system with dry friction,
J. Differential Equations  265 (2018) 5490--5519.

\bibitem{Chen4}
T. Chen, L. Huang, P. Yu,
Center condition and bifurcation of limit cycles for quadratic switching systems with a nilpotent equilibrium point,
J. Differential Equations 303 (2021) 326--368.

\bibitem{Chen2022}
T. Chen, J. Llibre,
Nilpotent center in a continuous piecewise quadratic polynomial Hamiltonian vector field,
Int. J. Bifurcation and Chaos 32 (2022) 2250116 (23 pages).

\bibitem{Cherkas1977}
L. A. Cherkas, Conditions for a Li\'enard equation to have a center,
Differential Equations 12(2) (1977) 201--206.

\bibitem{Christopher}
C. Christopher,
An algebraic approach to the classification of
centers in polynomial Li\'enard systems,
J. Appl. Math. Mech. 229 (1999) 329--329.

\bibitem{ChL1999}
C. Christopher, S. Lynch,
Small-amplitude limit cycle bifurcations for Li\'enard systems with quadratic or cubic
damping or restoring forces, Nonlinearity 12 (1999) 1099--1112.

\bibitem{CPG1999}
B. Coll, R. Prohens, A. Gasull,
The center problem for discontinuous Li\'enard differential equation,
Int. J. Bifurcation and Chaos 9 (1999) 1751--1761.

\bibitem{Colombo}
A. Colombo, P. Lamiani, L. Benadero, M. di Bernardo,
Two-parameter bifurcation analysis of the buck converter,
SIAM J. Appl. Dyn. Syst. 8 (2009) 1507--1522.

\bibitem{Dumortier2006}
F. Dumortier, J. Llibre, J. Art\'es,
Qualitative Theory of Planar Differential Systems,
Universitext, Springer-Verlag, New York, 2006.


\bibitem{Garcia}
I. Garc\'ia,
Cyclicity of some symmetric nilpotent centers,
J. Differential Equations 260 (2016) 5356--5377.

\bibitem{Gasull1989}
A. Gasull,
Differential equations that can be transformed into equations of Li\'enard type,
in ``$17^{\circ}$ Col\'oquio Brasileiro de Matem\'atica, 1989.''

\bibitem{Gasull}
A. Gasull,  J. Torregrosa,
Center-focus problem for discontinuous planar differential equations,
Int. J. Bifurcation and Chaos  13 (2003) 1755--1765.

\bibitem{Giacomini}
I. Garc\'ia, H. Giacomini, J. Gin\'e, J. Llibre,
Analytic nilpotent centers as limits of nondegenerate centers
revisited,
J. Math. Anal. Appl. 441 (2016) 893--899.


\bibitem{Han2013}
M. Han,
Bifurcation Theory of Limit Cycles, Science Press, Beijing, 2013.

\bibitem{HY2012}
M. Han, P. Yu,
Normal Forms, Melnikov Functions and Bifurcations of Limit Cycles,
Springer-Verlag, New York, 2012.

\bibitem{HK1991}
P. Hirschberg, E. Knobloch,
An unfolding of the Takens-Bogdanov singularity,
Quarterly Appl. Math. 49 (1991) 281--287.


\bibitem{JHYL2007}
J. Jiang, M. Han, P. Yu, S. Lynch,
Limit cycles in two types of symmetric Li\'enard systems,
Int. J. Bifurcation and Chaos 17 (2007) 2169--2174.

\bibitem{JH2009}
J. Jiang, M. Han,
Small-amplitude limit cycles of some Li\'enard-type systems,
Nonlinear Anal. 71 (2009) 6373--6377.


\bibitem{F.Li2021}
F. Li, H. Li, Y. Liu,
New Double Bifurcation of Nilpotent Focus,
Int. J. Bifurcation and Chaos 31 (2021) 2150053.

\bibitem{F.Li2}
F. Li, Y. Liu, Y. Liu, P. Yu,
Bi-center problem and bifurcation of limit cycles from nilpotent singular points in $Z_2$-equivariant cubic vector fields,
J. Differential Equations 265 (2018) 4965--4992.

\bibitem{LY2015}
F. Li, P. Yu, Y. Tian, Y. Liu,
Center and isochronous center conditions for switching systems associated with elementary singular points,
Commun. Nonlinear Sci. Numer. Simul. 28 (2015) 81--97.

\bibitem{LMP1977}
A. Lins, W. de Melo, C.C. Pugh,
On Li\'enard's equation, in: Geometry and Topology, in: Lecture Notes in Mathematics,
vol. 597, Springer, New York, 1977, pp. 335--357.

\bibitem{LT2015}
J. Llibre, M.A. Teixeira,
Limit cycles for m-piecewise discontinuous polynomial Li\'enard differential equations,
Z. Angew. Math. Phys. 66 (2015) 51--66.

\bibitem{LlibreV2022}
J. Llibre, C. Valls.
Global centers of the generalized polynomial Li\'enard differential systems,
J. Differential Equations 330 (2022) 66--80.


\bibitem{Liu2014}
Y. Liu, F. Li,
Double bifurcation of nilpotent focus,
Int. J. Bifurcation and Chaos 25 (2015) 1550036 (10 pages).

\bibitem{Liu2011}
Y. Liu, J. Li,
Bifurcations of limit cycles created by a multiple nilpotent critical
point of planar dynamical systems,
Int. J. Bifurcation and Chaos 21 (2011) 497--504.


\bibitem{Maesschalck2011}
P. De Maesschalck, F. Dumortier,
Classical Li\'enard equations of degree $n\geq6$ can have $[\frac{n-1}{2}]+2$ limit cycles,
J. Differential Equations 250 (2011) 2162--2176.

\bibitem{MM2014}
R.M. Martins, A.C. Mereu,
Limit cycles in discontinuous classical Li\'enard equations,
Nonlinear Anal., Real World Appl. 20 (2014) 67--73.


\bibitem{SHR2016}
L. Sheng, M. Han, V. Romanovsky,
On the number of limit cycles by perturbing a piecewise smooth Li\'enard model,
Int. J. Bifurcation and Chaos 26 (2016) 1650168.

\bibitem{Smale}
S. Smale,
Mathematical problems for the next century,
Math. Intelligencer 20 (1998) 7--15.

\bibitem{Strozyna}
E. Str\'{o}\.{z}zyna, H. {\rm \.{Z}o{\l}\c{a}dek},
The analytic and formal normal form for the nilpotent singularity,
J. Differential Equations 179 (2002) 479--537.

\bibitem{Tang}
S. Tang, J. Liang, Y. Xiao, R.A. Cheke,
Sliding bifurcations of Filippov two stage pest control models with economic
thresholds, SIAM J. Appl. Math. 72 (2012) 1061--1080.

\bibitem{TH2011}
Y. Tian, M. Han,
Hopf bifurcation for two types of Li\'enard systems,
J. Differential Equations 251 (2011) 834--859.

\bibitem{THX2019}
Y. Tian, M. Han, F. Xu,
Bifurcations of small limit cycles in Li\'enard systems with cubic restoring terms,
J. Differential Equations 267 (2019) 1561--1580.

\bibitem{TIAN}
Y. Tian,  P. Yu,
Center conditions in a switching Bautin system,
J. Differential Equations 259 (2015) 1203--1226.

\bibitem{TY2018}
Y. Tian, P. Yu,
Bifurcation of small limit cycles in cubic integrable systems using
higher-order analysis,
J. Differential Equations 264 (2018) 5950--5976.

\bibitem{Yang}
J. Yang, L. Zhao,
The cyclicity of period annuli for a class of cubic Hamiltonian systems
with nilpotent singular points,
J. Differential Equations 263 (2017) 5554--5581.

\bibitem{YHL2018}
P. Yu, M. Han, J. Li,
An improvement on the number of limit cycles bifurcating from a
non-degenerate center of homogeneous polynomial systems,
Int. J. Bifurcation and Chaos 28 (2018) 1850078 (31 pages).

\bibitem{YL2017}
P. Yu, F. Li,
Bifurcation of limit cycles in a cubic-order planar system around
a nilpotent critical point,
J. Math. Anal. Appl. 453(2) (2017) 645--667.

\bibitem{ZDHD1992}
Z. Zhang, T. Ding, W. Huang, Z. Dong,
Qualitative Theory of Differential Equations,
Transl. Math. Monogr., Amer.
Math. Soc., Providence, RI, 1992.

\end{thebibliography}
\end{document}